\documentclass[12pt]{article}
\usepackage{amssymb}
\usepackage{latexsym}
\usepackage{amsmath}

\newtheorem{definition}{Definition}[section]
\newtheorem{theorem}[definition]{Theorem}
\newtheorem{lemma}[definition]{Lemma}

\newtheorem{remark}[definition]{Remark}
\newtheorem{example}[definition]{Example}

\newtheorem{note}[definition]{Note}
\newtheorem{assumption}[definition]{Assumption}

\def\K{\mathbb K}

\def\Z{\mathbb Z}
\def\K{\mathbb K}

\begin{document}

\title{\bf Tridiagonal pairs and the $q$-tetrahedron algebra}
\author{Darren Funk-Neubauer\footnote{E-mail address: {\tt darren.funkneubauer@colostate-pueblo.edu } } \\
Department of Mathematics and Physics, Colorado State University-Pueblo \\ 2200 Bonforte Boulevard \\ Pueblo, CO 81001 USA}
\date{}

\maketitle

\begin{abstract}
The $q$-tetrahedron algebra $\boxtimes_q$ was recently introduced and has been studied in connection with tridiagonal pairs.  In this paper we further develop this connection.  Let $\K$ denote an algebraically closed field and let $q$ denote a nonzero scalar in $\K$ that is not a root of unity.  Let $V$ denote a vector space over $\K$ with finite positive dimension and let $A, A^*$ denote a tridiagonal pair on $V$.  Let $\lbrace \theta_i \rbrace _{i=0}^d$, (resp.  $\lbrace \theta^*_i \rbrace _{i=0}^d$) denote a standard ordering of the eigenvalues of $A$ (resp. $A^*$).  T. Ito and P. Terwilliger have shown that when $\theta_i = q^{2i-d}$ and $\theta^*_i = q^{d-2i}$ $(0 \leq i \leq d)$ there exists an irreducible $\boxtimes_q$-module structure on $V$ such that the $\boxtimes_q$ generators $x_{01}$, $x_{23}$ act as $A, A^*$ respectively.  In this paper we examine the case in which there exists a nonzero scalar $c$ in $\K$ such that $\theta_i = q^{2i-d}$ and $\theta^*_i = q^{2i-d} + cq^{d-2i}$ for $0 \leq i \leq d$.  In this case we associate to $A,A^*$ a polynomial $P$ in one variable and prove the following theorem as our main result.  \\
{\bf Theorem} The following are equivalent:
\begin{enumerate}
\item[{\rm (i)}]
There exists a $\boxtimes_q$-module structure on $V$ such that $x_{01}$ acts as $A$ and $x_{30} + cx_{23}$ acts as $A^*$, where $x_{01}, x_{30}, x_{23}$ are standard generators for $\boxtimes_q$.
\item[{\rm (ii)}]
$P( q^{2d-2} (q-q^{-1})^{-2} ) \neq 0$.
\end{enumerate}
Suppose (i),(ii) hold.  Then the $\boxtimes_q$-module structure on $V$ is unique and irreducible. \\ \\
{\bf AMS classification code:}  Primary: 17B37; Secondary: 15A21, 16W35, 17B65 \\ \\
{\bf Keywords:} tridiagonal pair, Leonard pair, tetrahedron algebra, $q$-tetrahedron algebra, quantum group, quantum affine algebra

\end{abstract}

\section{Tridiagonal pairs}

\noindent
We begin by recalling the notion of a {\it tridiagonal pair} \cite{ItoTanTer01}.  Let $\K$ denote a field and let $V$ denote a vector space over $\K$ with finite positive dimension.  Let $A:V \rightarrow V$ denote a linear transformation and let $W$ denote a subspace of $V$.  We say $W$ is an {\it eigenspace} of $A$ whenever $W \neq 0$ and there exists $\theta \in \K$ such that
\begin{eqnarray*}
W= \{ v \in V \, | \, Av=\theta v \, \}.
\end{eqnarray*}
In this case, we call $\theta$ an {\it eigenvalue} of $A$.  We say $A$ is {\it diagonalizable} whenever $V$ is spanned by the eigenspaces of $A$.

\begin{definition}
\rm
\cite[Definition 1.1]{ItoTanTer01}
\label{def:tdp}
Let $V$ denote a vector space over $\K$ with finite positive dimension.  By a {\it tridiagonal pair} on $V$,
we mean an ordered pair of linear transformations $A:V\to V$ and $A^*:V \to V$ that satisfy the following four conditions:
\begin{enumerate}
\item Each of $A,A^*$ is diagonalizable.
\item There exists an ordering $\lbrace V_i \rbrace _{i=0}^d$ of the eigenspaces of $A$ such that
\begin{equation}
A^* V_i \subseteq V_{i-1} + V_i+ V_{i+1} \qquad \qquad (0 \leq i \leq d),
\label{eq:t1}
\end{equation}
where $V_{-1} = 0$, $V_{d+1}= 0$.
\item There exists an ordering $\lbrace V^*_i \rbrace _{i=0}^{\delta}$ of the eigenspaces of $A^*$ such that
\begin{equation}
A V^*_i \subseteq V^*_{i-1} + V^*_i+ V^*_{i+1} \qquad \qquad (0 \leq i \leq \delta),
\label{eq:t2}
\end{equation}
where $V^*_{-1} = 0$, $V^*_{\delta+1}= 0$.
\item There does not exist a subspace $W$ of $V$ such  that $AW\subseteq W$, $A^*W\subseteq W$, $W\not=0$, $W\not=V$.
\end{enumerate}
\end{definition}

\begin{note}
\rm
According to a common notational convention $A^*$ denotes the conjugate-transpose of $A$.  We are not using this convention.  For a tridiagonal pair $A,A^*$ the linear transformations $A$ and $A^*$ are arbitrary subject to (i)--(iv) above.
\end{note}

\noindent
Tridiagonal pairs originally arose in algebraic combinatorics through the study of a combinatorial object called a P- and Q-polynomial association scheme \cite{ItoTanTer01}.  Since then they have appeared in many other areas of mathematics.  For instance, examples of tridiagonal pairs appear in representation theory \cite{Al-Najjar05, Benkart04, Bowman07, Funk-Neubauer07, ItoTer073, ItoTer075}, the study of orthogonal polynomials and special functions \cite{Ter011, Ter041, Ter042}, the theory of partially ordered sets \cite{Ter90, Ter03}, and statistical mechanics \cite{Baseilhac05, HarTer07, Ter012}.  The tridiagonal pairs for which the $V_i, V^*_i$ all have dimension 1 are called {\it Leonard pairs}.  The Leonard pairs are classified and correspond to a family of orthogonal polynomials consisting of the $q$-Racah polynomials and related polynomials in the Askey scheme \cite{Ter011, Ter041}.  Currently there is no classification of tridiagonal pairs; this paper is largely motivated by the search for this classification.  For further information on tridiagonal pairs and Leonard pairs see \cite{Al-Najjar04, Al-Najjarinpress, Baseilhac06, Curtin071, Curtin072, Hartwig05, ItoTer04, ItoTer074, Nomura051, Nomura052, Nomura053, Nomurainpress2, Nomura061,  Nomura062, Nomura063, Nomura071, Nomura072, Nomura073, Nomura074, Nomurainpress1, Nomura081, Nomurainpress3, Ter02, Ter051, Ter052, Ter053, TerVid04, Vidarinpress, Vid07}.

\medskip

\noindent
We now recall a few basic facts about tridiagonal pairs.  Let $A,A^*$ denote a tridiagonal pair on $V$ and let $d, \delta$ be as in Definition \ref{def:tdp}(ii), (iii).  By \cite[Lemma~4.5]{ItoTanTer01} we have $d=\delta$; we call this common value the {\it diameter} of $A,A^*$.  An ordering of the eigenspaces of $A$ (resp. $A^*$) will be called {\it standard} whenever it satisfies (\ref{eq:t1}) (resp. (\ref{eq:t2})).  We comment on the uniqueness of the standard ordering.  Let $\lbrace V_i \rbrace _{i=0}^d$ denote a standard ordering of the eigenspaces of $A$.  Then the ordering $\lbrace V_{d-i} \rbrace _{i=0}^d$ is standard and no other ordering is standard.  A similar result holds for the eigenspaces of $A^*$.  An ordering of the eigenvalues of $A$ (resp. $A^*$) will be called {\it standard} whenever the corresponding ordering of the eigenspaces of $A$ (resp. $A^*$) is standard.  Let  $\lbrace \theta_i \rbrace _{i=0}^d$ (resp.  $\lbrace \theta^*_i \rbrace _{i=0}^d$) denote a standard ordering of the eigenvalues of $A$ (resp. $A^*$).  The $\theta_i, \theta^*_i$ both satisfy a three term recurrence relation that has been solved in closed form \cite[Theorem 11.1, Theorem 11.2]{ItoTanTer01}.  The following two special cases will be of interest to us.  For $0 \neq q \in \K$ we call $A,A^*$ {\it $q$-geometric} whenever $\theta_i=q^{2i-d}$ and $\theta^*_i=q^{d-2i}$ for $0 \leq i \leq d$.  For a study of this case see \cite{Al-Najjar05, Al-Najjarinpress, ItoTer073, ItoTer075}.  For $0 \neq q \in \K$ we call $A,A^*$ {\it $q$-mixed} whenever there exists $0 \neq c \in \K$ such that $\theta_i=q^{2i-d}$ and $\theta^*_i=q^{2i-d} + c \, q^{d-2i}$ for $0 \leq i \leq d$.  The main result of this paper concerns the $q$-mixed tridiagonal pairs.  These two cases are of interest because of their connection to the $q$-tetrahedron algebra $\boxtimes_q$.  We discuss this connection in the next section.

\section{The $q$-tetrahedron algebra}

\noindent
The $q$-tetrahedron algebra $\boxtimes_q$ was introduced in \cite{ItoTer072} as part of the continuing investigation of tridiagonal pairs.  It is closely related to a number of well known algebras including the quantum group $U_q(\mathfrak{sl}_2)$ \cite[Proposition 7.4]{ItoTer072}, the $U_q(\mathfrak{sl}_2)$ loop algebra \cite[Proposition 8.3]{ItoTer072}, and positive part of $U_q(\widehat{\mathfrak{sl}}_2)$ \cite[Proposition 9.4]{ItoTer072}.  The finite dimensional irreducible $\boxtimes_q$-modules are described in \cite{ItoTer072}.  For further information on $\boxtimes_q$ see \cite{ItoTer071}.  We note that $\boxtimes_q$ is a $q$-analogue of the tetrahedron algebra $\boxtimes$ \cite{Benkart07, Elduque07, Hartwig07, HarTer07, ItoTerinpress4}.

\medskip

\noindent
We now fix some notation and then recall the definition of $\boxtimes_q$.

\medskip

\noindent
Throughout the rest of this paper $\K$ will denote an algebraically closed field.  We fix a nonzero scalar $q \in \K$ that is not a root of $1$.  For an integer $n \geq 0$ we define
\begin{eqnarray}
\label{def:brackn}
\lbrack n \rbrack = \frac{q^n-q^{-n}}{q-q^{-1}}.
\end{eqnarray}

\noindent
We let $\Z_4 = \Z/4\Z$ denote the cyclic group of order 4.

\begin{definition}
\rm
\cite[Definition 6.1]{ItoTer072}
\label{def:qtet}
Let $\boxtimes_q$ denote the unital associative $\K$-algebra that has generators\begin{eqnarray*}
\lbrace x_{ij}\;|\; i,j \in \Z_4,\;j-i=1 \;\mbox{or} \;j-i=2\rbrace
\end{eqnarray*}
and the following relations:
\begin{enumerate}
\item For $i,j\in \Z_4$ such that $j-i=2$,
\begin{eqnarray*}
x_{ij}x_{ji} = 1.
\label{eq:qrel0}
\end{eqnarray*}
\item For $h,i,j\in \Z_4$ such that the pair $(i-h,j-i)$ is one of $(1,1), (1,2), (2,1)$,
\begin{eqnarray}
\frac{qx_{ij}x_{hi}-q^{-1}x_{hi}x_{ij}}{q-q^{-1}}=1.
\label{eq:qrel1}
\end{eqnarray}
\item For $h,i,j,k\in \Z_4$ such that $i-h=j-i=k-j=1$,
\begin{eqnarray}
\label{eq:qserre}
x_{hi}^3x_{jk} - \lbrack 3 \rbrack x_{hi}^2x_{jk}x_{hi} + \lbrack 3 \rbrack x_{hi}x_{jk}x_{hi}^2 - x_{jk}x_{hi}^3=0.
\end{eqnarray}
\end{enumerate}
We call $\boxtimes_q$ the {\it $q$-tetrahedron algebra}.  We refer to the $x_{ij}$ as the {\it standard generators} for $\boxtimes_q$.
\end{definition}

\begin{remark}
\rm
The equations  (\ref{eq:qserre}) are called the {\it cubic $q$-Serre relations}.\end{remark}

\noindent
We now recall a few basic facts about $\boxtimes_q$-modules.  Let $V$ denote a finite dimensional irreducible $\boxtimes_q$-module.  By \cite[Theorem 12.3]{ItoTer072} each generator $x_{ij}$ of $\boxtimes_q$ is diagonalizable on $V$.  Moreover, there exist an integer $d\geq 0$ and a scalar $\varepsilon \in \lbrace 1,-1\rbrace$ such that for each generator $x_{ij}$ the set of distinct eigenvalues of $x_{ij}$ on $V$ is $\lbrace \varepsilon q^{2n-d} \,|\, 0 \leq n \leq d \rbrace$.  We call $\varepsilon$ the {\it type} of $V$.

\medskip

\noindent
We now discuss the connection between $q$-geometric tridiagonal pairs and finite dimensional irreducible $\boxtimes_q$-modules.  Let $V$ denote a vector space over $\K$ with finite positive dimension.  Let $A,A^*$ denote a $q$-geometric tridiagonal pair on $V$.  Then there exists an irreducible $\boxtimes_q$-module structure on $V$ of type $1$ such that $A$ acts as $x_{01}$ and $A^*$ acts as $x_{23}$.  Conversely, let $V$ denote a finite dimensional irreducible $\boxtimes_q$-module of type $1$.  Then the generators $x_{01},x_{23}$ act on $V$ as a $q$-geometric tridiagonal pair \cite[Theorem 34.14]{Ter042}.

\medskip

\noindent
Inspired by the above result we consider the connection between $q$-mixed tridiagonal pairs and finite dimensional irreducible $\boxtimes_q$-modules.  In the next section we give a detailed description of the situation we wish to consider.

\section{Assumptions and motivation}

\noindent
Throughout the rest of the paper we will be concerned with the following situation.

\begin{assumption}
\label{assum}
\rm
Let $V$ denote a vector space over $\K$ with finite positive dimension and let $A,A^*$ denote a tridiagonal pair on $V$.  Let $\lbrace V_i \rbrace _{i=0}^d$ (resp. $\lbrace V^*_i \rbrace _{i=0}^d$) denote a standard ordering of the eigenspaces of $A$ (resp. $A^*$).  For $0 \leq i \leq d$, let $\theta_i$ (resp. $\theta^*_i$) denote the eigenvalue of $A$ (resp. $A^*$) associated with $V_i$ (resp. $V^*_i$).  We assume there exists a nonzero $c \in \K$ such that $\theta_i=q^{2i-d}$ and $\theta^*_i=q^{2i-d} + c \, q^{d-2i}$ for $0 \leq i \leq d$.
\end{assumption}

\noindent
To motivate our main result we make some comments.

\begin{lemma}
\label{thm:AA^*rel}
{\rm \cite[Theorem 10.1, Theorem 11.1]{ItoTanTer01}}
With reference to Assumption \ref{assum} we have
\begin{enumerate}
\item $A^3A^* - \lbrack 3 \rbrack A^2A^*A + \lbrack 3 \rbrack AA^*A^2 - A^*A^3=0$,
\item $A^{* \, 3}A - \lbrack 3 \rbrack A^{* \, 2}AA^* + \lbrack 3 \rbrack A^*AA^{* \, 2} - AA^{* \, 3} + c(q^2-q^{-2})^2(A^*A-AA^*)=0$.
\end{enumerate}
\end{lemma}

\noindent
  Using (\ref{eq:qrel1}) and (\ref{eq:qserre}) it can be shown that for $0 \neq c \in \K$ the elements $x_{01}$ and $x_{30}+cx_{23}$ of $\boxtimes_q$ satisfy the relations in Lemma \ref{thm:AA^*rel}.  Given this, it is natural to ask the following question.  With reference to Assumption \ref{assum} when does there exist an irreducible $\boxtimes_q$-module structure on $V$ of type $1$ such that $A$ acts as $x_{01}$ and $A^*$ acts as $x_{30} + c x_{23}$?  In this paper we answer this question; our main result is Theorem \ref{thm:main}.  In the next section we establish some notation needed to state our main result.

\section{A split decomposition and its raising/lowering maps}

\noindent
We now recall the notion of a split decomposition of a tridiagonal pair and its corresponding raising and lowering maps.

\begin{definition}
\label{def:decomp}
\rm
Let $V$ denote a vector space over $\K$ with finite positive dimension.  By a \\{\it decomposition} of $V$ we mean a sequence  $\lbrace U_i \rbrace _{i=0}^d$ consisting of nonzero subspaces of $V$ such that $V=\sum_{i=0}^{d} U_i$ (direct sum).  For notational convenience we set $U_{-1}:=0, U_{d+1}:=0$.
\end{definition}

\noindent
Referring to Assumption \ref{assum} the sequences $\lbrace V_i \rbrace _{i=0}^d$ and  $\lbrace V^*_i \rbrace _{i=0}^d$ are both decompositions of $V$.  We now mention another decomposition of interest.

\begin{lemma}
\label{def:U_i}
{\rm \cite[Theorem~4.6]{ItoTanTer01}}
With reference to Assumption \ref{assum}, for $0 \leq i \leq d$ define
\begin{eqnarray*}
U_i=(V^*_0 + \cdots + V^*_i) \cap (V_i + \cdots + V_d).
\end{eqnarray*}
Then $\lbrace U_i \rbrace _{i=0}^d$ is a decomposition of $V$.  Moreover, for $0 \leq i \leq d$
\begin{align}
\label{AA^*U1}
(A^*-\theta^*_iI)U_i &\subseteq U_{i-1}, & (A-\theta_iI)U_i &\subseteq U_{i+1}, \\
\label{AA^*U2}
U_0 + \cdots + U_i &= V^*_0 + \cdots + V^*_i,  & U_i + \cdots + U_d &= V_i + \cdots + V_d.
\end{align}
We call $\lbrace U_i \rbrace _{i=0}^d$ the {\rm split decomposition} of $V$ corresponding to the given orderings $\lbrace V_i \rbrace _{i=0}^d$, $\lbrace V^*_i \rbrace _{i=0}^d$.
\end{lemma}

\begin{definition}
\label{def:Fi}
{\rm \cite[Definition~5.2]{ItoTanTer01}}
\rm
With reference to Assumption \ref{assum} and Lemma \ref{def:U_i} we define the following.  For $0 \leq i \leq d$ we define a linear transformation $F_i: V \rightarrow V$ by
\begin{center}
$(F_i - I)U_i=0$, \\
$F_iU_j=0 \,\,\,\,  \mbox{if} \,\,\,\, j \neq i, \qquad (0 \leq j \leq d)$.
\end{center}
In other words, $F_i$ is the projection map from $V$ onto $U_i$.  We call $F_i$ the {\it $i$th projection map corresponding to $\lbrace U_i \rbrace _{i=0}^d$}.
\end{definition}

\begin{definition}
\label{def:RL}
{\rm \cite[Definition~6.1]{ItoTanTer01}}
\rm
With reference to Assumption \ref{assum} and Definition \ref{def:Fi} we define
\begin{eqnarray}
\label{RL1}
R = A - \sum_{h=0}^{d}\theta_h F_h, \qquad \qquad L = A^* - \sum_{h=0}^{d}\theta^*_h F_h.
\end{eqnarray}
Combining (\ref{AA^*U1}) and (\ref{RL1}) we have $RU_i \subseteq U_{i+1}$ and $LU_i \subseteq U_{i-1}$ for $0 \leq i \leq d$.  We call $R$ (resp. $L$) the {\it raising (resp. lowering) map corresponding to $\lbrace U_i \rbrace _{i=0}^d$}.
\end{definition}

\section{The main theorem}

\noindent
In this section we state our main result.  We begin with a few comments.

\begin{lemma}
\label{thm:dim1U0}
{\rm \cite[Theorem~1.3]{Nomurainpress2}}
With reference to Assumption \ref{assum} and Lemma \ref{def:U_i} we have \rm{dim}$(U_0)=1$.
\end{lemma}

\begin{definition}
\label{def:zeta}
\rm
With reference to Definition \ref{def:RL} and Lemma \ref{thm:dim1U0} we find that for $0 \leq i \leq d$ $U_0$ is contained in an eigenspace for $L^iR^i$; let $\zeta_i$ denote the corresponding eigenvalue.
\end{definition}

\noindent
With reference to (\ref{def:brackn}) for an integer $n \geq 0$ we define
\begin{eqnarray}
\label{def:brackfacn}
\lbrack n \rbrack ! = \lbrack n \rbrack \lbrack n-1 \rbrack \cdots \lbrack 1 \rbrack.
\end{eqnarray}
We interpret $\lbrack 0 \rbrack ! = 1$.

\begin{definition}
\label{def:PV}
\rm
With reference to Assumption \ref{assum} and Definition \ref{def:zeta} we define a polynomial $P \in \K\lbrack \lambda \rbrack$ ($\lambda$ indeterminate) by
\begin{eqnarray*}
P = \sum_{i=0}^{d} \frac{q^{i(1-i)} \, \zeta_i \, \lambda^i}{\lbrack i \rbrack ! \, ^2}.
\end{eqnarray*}
\end{definition}

\noindent
We now state our main result.

\begin{theorem}
\label{thm:main}
With reference to Assumption \ref{assum} the following are equivalent:
\begin{enumerate}
\item[{\rm (i)}]
There exists a $\boxtimes_q$-module structure on $V$ such that $x_{01}$ acts as $A$ and $x_{30} + cx_{23}$ acts as $A^*$.
\item[{\rm (ii)}]
$P( q^{2d-2} (q-q^{-1})^{-2} ) \neq 0$ where $P$ is from Definition \ref{def:PV}.
\end{enumerate}
Suppose (i),(ii) hold.  Then the $\boxtimes_q$-module structure on $V$ is unique, irreducible, and has type $1$.
\end{theorem}

\section{An outline of the proof of Theorem \ref{thm:main}}

\noindent
Our proof of Theorem \ref{thm:main} will consume the remainder of the paper from Section 7 to Section 18.  Here we sketch an overview of the argument.

\medskip

\noindent
We adopt Assumption \ref{assum}.  The main idea used in proving Theorem \ref{thm:main} is the following.  We modify the linear transformation $A^*:V \rightarrow V$ to produce a new linear transformation $\widetilde{A}^*:V \rightarrow V$ and we show that $A,\widetilde{A}^*$ is a $q$-geometric tridiagonal pair on $V$ if and only if $P( q^{2d-2} (q-q^{-1})^{-2} ) \neq 0$.  Then we apply \cite[Theorem 2.7]{ItoTer075} and \cite[Theorem 10.4]{ItoTer072} to $A, \widetilde{A}^*$ to produce the $\boxtimes_q$-module structure on $V$ as in Theorem \ref{thm:main}.  The plan for the paper is as follows.  In Sections 7 and 8 we present some lemmas and definitions which will be used as tools throughout the remainder of the paper.  In Section 9 we define the linear transformation $\widetilde{A}^*:V \rightarrow V$.  We show that $\widetilde{A}^*$ is diagonalizable on $V$ and the set of distinct eigenvalues of $\widetilde{A}^*$ on $V$ is $\lbrace q^{d-2i} \,|\, 0 \leq i \leq d \rbrace$.  In Section 10 we show that $A$ and $\widetilde{A}^*$ satisfy Definition \ref{def:tdp}(ii),(iii).  Sections 11 through 17 are devoted to showing that $A$ and $\widetilde{A}^*$ satisfy Definition \ref{def:tdp}(iv) if and only if $P( q^{2d-2} (q-q^{-1})^{-2} ) \neq 0$.  We note that the arguments given in Sections 11 through 17 are a modification of the arguments from \cite[Sections 7--12]{ItoTer075}.  In Section 18 we show how to use \cite[Theorem 2.7]{ItoTer075} and \cite[Theorem 10.4]{ItoTer072} applied to the $q$-geometric tridiagonal pair $A,\widetilde{A}^*$ to produce the $\boxtimes_q$-module structure on $V$ as in Theorem \ref{thm:main}.

\section{Some more raising/lowering maps}

\noindent
We now present another split decomposition for the tridiagonal pair $A,A^*$ and its corresponding raising and lowering maps.

\begin{definition}
\label{def:W_iandfriends}
\rm
With reference to Assumption \ref{assum} and Lemma \ref{def:U_i} let  $\lbrace W_i \rbrace _{i=0}^d$ denote the split decomposition of $V$ corresponding to the orderings $\lbrace V_{d-i} \rbrace _{i=0}^d$, $\lbrace V^*_i \rbrace _{i=0}^d$.  With reference to Definition \ref{def:Fi} for $0 \leq i \leq d$, let $G_i$ denote the $i$th projection map corresponding to $\lbrace W_i \rbrace _{i=0}^d$.  With reference to Definition \ref{def:RL} let $r$ (resp. $l$) denote the raising (resp. lowering) map corresponding to $\lbrace W_i \rbrace _{i=0}^d$.
\end{definition}

\noindent
We make the following three remarks in order to emphasize the similarities and differences between the two split decompositions $\lbrace U_i \rbrace _{i=0}^d$ and $\lbrace W_i \rbrace _{i=0}^d$.

\begin{remark}
\label{W_iformulas}
\rm
With reference to Definition \ref{def:W_iandfriends} we emphasize the following.  For $0 \leq i \leq d$ we have
\begin{eqnarray*}
W_i=(V^*_0 + \cdots + V^*_i) \cap (V_0 + \cdots + V_{d-i}).
\end{eqnarray*}
Moreover, for $0 \leq i \leq d$
\begin{align}
\label{AA^*W1}
(A^*-\theta^*_iI)W_i &\subseteq W_{i-1}, & (A-\theta_{d-i}I)W_i &\subseteq W_{i+1}, \\
\label{AA^*W2}
W_0 + \cdots + W_i &= V^*_0 + \cdots + V^*_i,  & W_i + \cdots + W_d &= V_0 + \cdots + V_{d-i}.
\end{align}
\end{remark}

\begin{remark}
\label{G_iformulas}
\rm
With reference to Definition \ref{def:W_iandfriends} we emphasize that for $0 \leq i \leq d$
\begin{center}
$(G_i - I)W_i = 0$, \\
$G_iW_j = 0 \,\,\,\,  \mbox{if} \,\,\,\, j \neq i, \qquad (0 \leq j \leq d)$.
\end{center}
\end{remark}

\begin{remark}
\label{rlformulas}
\rm
With reference to Definition \ref{def:W_iandfriends} we emphasize that
\begin{eqnarray}
\label{rl}
r = A - \sum_{h=0}^{d}\theta_{d-h} G_h, \qquad \qquad l = A^* - \sum_{h=0}^{d}\theta^*_h G_h.
\end{eqnarray}
Moreover, for $0 \leq i \leq d$, $rW_i \subseteq W_{i+1}$ and $lW_i \subseteq W_{i-1}$.
\end{remark}

\section{Some linear algebra}

\noindent
In this section we state some linear algebraic results that will be useful throughout the paper.

\medskip

\noindent
We use the following notation.  Let $V$ denote a finite dimensional vector space over $\K$ and let $X:V \rightarrow V$ denote a linear transformation.  For $\theta \in \K$ we define \begin{eqnarray*}
V_{X}(\theta) = \{ v \in V \, | \, X v = \theta v \, \}.
\end{eqnarray*}
Observe that $\theta$ is an eigenvalue of $X$ if and only if $V_{X}(\theta) \neq 0$, and in this case $V_{X}(\theta)$ is the corresponding eigenspace.

\begin{lemma}
\label{thm:bidiag1}
{\rm \cite[Lemma~11.2]{ItoTer072}}
Let $V$ denote a vector space over $\K$ with finite positive dimension.  Let $X: V \rightarrow V$ and $Y: V \rightarrow V$ denote linear transformations.  Then for all nonzero $\theta \in \K$ the following are equivalent:
\begin{enumerate}
\item[{\rm (i)}]
The expression $qXY-q^{-1}YX-(q-q^{-1})I$ vanishes on $V_X(\theta)$.
\item[{\rm (ii)}]
$(Y -\theta^{-1} I)V_X(\theta) \subseteq V_X(q^{-2}\theta)$.
\end{enumerate}
\end{lemma}

\begin{lemma}
\label{thm:bidiag2}
{\rm \cite[Lemma~11.3]{ItoTer072}}
Let $V$ denote a vector space over $\K$ with finite positive dimension. Let $X: V \rightarrow V$ and $Y: V \rightarrow V$ denote linear transformations.  Then for all nonzero $\theta \in \K$ the following are equivalent:
\begin{enumerate}
\item[{\rm (i)}]
The expression $qXY-q^{-1}YX-(q-q^{-1})I$ vanishes on $V_Y(\theta)$.
\item[{\rm (ii)}]
$(X -\theta^{-1}I)V_Y(\theta) \subseteq V_Y(q^{2}\theta)$.
\end{enumerate}
\end{lemma}

\begin{lemma}
\label{thm:bidiag3}
Let $V$ denote a vector space over $\K$ with finite positive dimension.  Let $X: V \rightarrow V$ and $Y: V \rightarrow V$ denote linear transformations.  Fix a nonzero $c \in \K$.  Then for all nonzero $\theta \in \K$ the following are equivalent:
\begin{enumerate}
\item[{\rm (i)}]
The expression $qXY-q^{-1}YX-(q-q^{-1})(X^2+c \, I)$ vanishes on $V_X(\theta)$.
\item[{\rm (ii)}]
$(Y -\theta I-c \, \theta^{-1} I)V_X(\theta) \subseteq V_X(q^{-2}\theta)$.
\end{enumerate}
\end{lemma}

\noindent
{\it Proof:} For $v \in V_X(\theta)$ we have
\begin{eqnarray*}
(qXY-q^{-1}YX-(q-q^{-1})(X^2+c \, I))v=q(X-q^{-2}\theta I)(Y-\theta I-c \, \theta^{-1} I)v
\end{eqnarray*}
and the result follows.
\hfill $\Box $ \\

\begin{lemma}
\label{thm:qweylsum}
{\rm \cite[Lemma~11.4]{ItoTer072}}
Let $V$ denote a vector space over $\K$ with finite positive dimension.  Let $X:V\rightarrow V$ and $Y:V \rightarrow V$ denote linear transformations such that
\begin{eqnarray*}
\frac{qXY-q^{-1}YX}{q-q^{-1}}=I.
\end{eqnarray*}
Then for all nonzero $\theta \in \K$,
\begin{eqnarray}
\sum_{n=0}^\infty V_X(q^{-2n}\theta)=\sum_{n=0}^\infty V_Y(q^{2n}\theta^{-1}).
\label{qweylsum1}
\end{eqnarray}
\end{lemma}

\section{The linear transformations $B$ and $\widetilde{A}^*$}

\begin{definition}
\label{def:B}
\rm
With reference to Assumption \ref{assum} and Definition \ref{def:W_iandfriends} let $B:V \rightarrow V$ denote the linear transformation such that for $0 \leq i \leq d$, $W_i$ is an eigenspace of $B$ with eigenvalue $q^{2i-d}$.
\end{definition}

\begin{lemma}
\label{thm:AA^*B}
With reference to Assumption \ref{assum} and Definition \ref{def:B} we have
\begin{eqnarray}
\frac{qAB-q^{-1}BA}{q-q^{-1}} &=& I,
\label{AA^*B1} \\
\frac{qBA^*-q^{-1}A^*B}{q-q^{-1}} &=& B^2 + c \, I.
\label{AA^*B2}
\end{eqnarray}
\end{lemma}

\noindent
{\it Proof:} Recall that $\lbrace W_i \rbrace _{i=0}^d$ is a decomposition of $V$.  By (\ref{AA^*W1}) $(A-q^{d-2i}I)W_i \subseteq W_{i+1}$ for $0 \leq i \leq d$.  Using this and Lemma \ref{thm:bidiag2} we obtain (\ref{AA^*B1}).  By (\ref{AA^*W1}) $(A^*-q^{2i-d}-cq^{d-2i}I)W_i \subseteq W_{i-1}$ for $0 \leq i \leq d$.  Using this and Lemma \ref{thm:bidiag3} we obtain (\ref{AA^*B2}).
\hfill $\Box$ \\

\begin{definition}
\label{def:A^*tilde}
\rm
With reference to Assumption \ref{assum} and Definition \ref{def:B} let $\widetilde{A}^*: V \rightarrow V$ denote the following linear transformation:
\begin{eqnarray*}
\widetilde{A}^*=c^{-1}(A^* - B).
\end{eqnarray*}
\end{definition}

\begin{lemma}
\label{thm:A^*tildeW}
With reference to Definition \ref{def:W_iandfriends} and Definition \ref{def:A^*tilde} we have
\begin{eqnarray*}
(\widetilde{A}^*-q^{d-2i}I)W_i \subseteq W_{i-1} \qquad \qquad (0 \leq i \leq d)
\end{eqnarray*}
\end{lemma}

\noindent
{\it Proof:}  Let $i$ be given.  Recall that $W_i$ is an eigenspace for $B$ with eigenvalue $q^{2i-d}$.  We have
\[ \begin{array}{cccc}
(\widetilde{A}^*-q^{d-2i}I) W_i &=&
c^{-1}(A^*-q^{2i-d}I-c \, q^{d-2i}I) W_i \hfill&  \,\,\,\,\,\, (\mbox{by Definition \ref{def:A^*tilde}})
\\
&\subseteq&
W_{i-1} \hfill&\hskip 1 truein  (\mbox{by (\ref{AA^*W1})}
).
\end{array}\]
\hfill $\Box$ \\

\begin{lemma}
\label{thm:A^*tildeBrel}
With reference to Definition \ref{def:B} and Definition \ref{def:A^*tilde} we have
\begin{eqnarray}
\label{A^*tildeBrel1}
\frac{qB \widetilde{A}^*-q^{-1} \widetilde{A}^* B}{q-q^{-1}} = I.
\end{eqnarray}
\end{lemma}

\noindent
{\it Proof:} Immediate from Lemma \ref{thm:bidiag1} and Lemma \ref{thm:A^*tildeW}.
\hfill $\Box$ \\

\begin{lemma}
\label{thm:A^*tildediag}
With reference to Definition \ref{def:A^*tilde} the following holds.  $\widetilde{A}^*$ is diagonalizable with eigenvalues $\lbrace q^{d-2i} \rbrace _{i=0}^d$.  Moreover, for $0 \leq i \leq d$, the dimension of the eigenspace of $\widetilde{A}^*$ associated with $q^{d-2i}$ is equal to the dimension of $W_i$.
\end{lemma}

\noindent
{\it Proof:}  We start by displaying the eigenvalues of $\widetilde{A}^*$. Notice that the scalars $q^{d-2i}$ ($0 \leq i \leq d$) are distinct since $q$ is not a root of unity.  Using Lemma \ref{thm:A^*tildeW} we see that, with respect to an appropriate basis for $V$, $\widetilde{A}^*$ is represented by a upper triangular matrix that has diagonal entries $q^{d}, q^{d-2}, \ldots, q^{-d}$, with $q^{d-2i}$ appearing dim($W_i$) times for $0 \leq i \leq d$.  Hence for $0 \leq i \leq d$ $q^{d-2i}$ is a root of the characteristic polynomial of $\widetilde{A}^*$ with multiplicity dim$(W_i)$.  It remains to show that $\widetilde{A}^*$ is diagonalizable.  To do this we show that the minimal polynomial of $\widetilde{A}^*$ has distinct roots.  Recall that $\lbrace W_i \rbrace _{i=0}^d$ is a decomposition of $V$.  Using Lemma \ref{thm:A^*tildeW} we find that $\prod_{i=0}^{d}(\widetilde{A}^*-q^{d-2i}I)V=0$.  By this and since $q^{d-2i}$ ($0 \leq i \leq d$) are distinct we see that the minimal polynomial of $\widetilde{A}^*$ has distinct roots.  We conclude that $\widetilde{A}^*$ is diagonalizable and the result follows.
\hfill $\Box$ \\

\begin{definition}
\label{def:V^*_itilde}
\rm
With reference to Definition \ref{def:A^*tilde} and Lemma \ref{thm:A^*tildediag}, for $0 \leq i \leq d$ we let $\widetilde{V}^*_i$ denote the eigenspace for $\widetilde{A}^*$ with eigenvalue $q^{d-2i}$.  For notational convenience we set $\widetilde{V}^*_{-1}:=0,\widetilde{V}^*_{d+1}:=0$.  We observe that  $\lbrace \widetilde{V}^*_i \rbrace _{i=0}^d$ is a decomposition of $V$.
\end{definition}

\section{The linear transformations $A$, $\widetilde{A}^*$ satisfy the cubic $q$-Serre relations}

\begin{lemma}
\label{thm:AA^*tilderel}
With reference to Assumption \ref{assum} and Definition \ref{def:A^*tilde} we have
\begin{enumerate}
\item $A^3 \widetilde{A}^*-\lbrack 3 \rbrack A^{2} \widetilde{A}^*A+\lbrack 3 \rbrack A \widetilde{A}^*A^2-\widetilde{A}^*A^3=0$,
\item $\widetilde{A}^{*3}A-\lbrack 3 \rbrack \widetilde{A}^{*2}A \widetilde{A}^*+\lbrack 3 \rbrack \widetilde{A}^*A \widetilde{A}^{*2}-A \widetilde{A}^{*3}=0$.
\end{enumerate}
\end{lemma}

\noindent
{\it Proof:}  By Definition \ref{def:A^*tilde} we have $A^*=c \widetilde{A}^* + B$.  Substitute this into Lemma \ref{thm:AA^*rel}(i),(ii) and simplify the result using (\ref{AA^*B1}) and (\ref{A^*tildeBrel1}).
\hfill $\Box$ \\

\begin{lemma}
\label{thm:tridiagAA^*tilde}
With reference to Assumption \ref{assum}, Definition \ref{def:A^*tilde}, and Definition \ref{def:V^*_itilde} we have
\begin{enumerate}
\item $\widetilde{A}^* V_i \subseteq V_{i-1}+V_i+V_{i+1}$, \qquad $0 \leq i \leq d$,
\item  $A \widetilde{V}^*_i \subseteq \widetilde{V}^*_{i-1}+ \widetilde{V}^*_i + \widetilde{V}^*_{i+1}$, \qquad $0 \leq i \leq d$.
\end{enumerate}
\end{lemma}

\noindent
{\it Proof:} Immediate from Lemma \ref{thm:AA^*tilderel} and \cite[Lemma~11.1]{ItoTer072}.
\hfill $\Box$ \\

\begin{remark}
\rm
Recall that in order to prove Theorem \ref{thm:main} we need to show that $A,\widetilde{A}^*$ is a $q$-geometric tridiagonal pair on $V$ if and only if $P( q^{2d-2} (q-q^{-1})^{-2} ) \neq 0$ (see Theorem \ref{thm:AA^*tildetrip}).  Combining Assumption \ref{assum}, Lemma \ref{thm:A^*tildediag}, and Lemma \ref{thm:tridiagAA^*tilde} we have that $A, \widetilde{A}^*$ satisfy Definition \ref{def:tdp}(i),(ii),(iii).  Sections 11 through 17 are devoted to showing that $A, \widetilde{A}^*$ satisfy Definition \ref{def:tdp}(iv) if and only if $P( q^{2d-2} (q-q^{-1})^{-2} ) \neq 0$ (see Theorem \ref{thm:PVnot0}).
\end{remark}

\section{The linear transformation $K$}

\begin{definition}
\label{def:K}
\rm
With reference to Assumption \ref{assum} and Lemma \ref{def:U_i} let $K:V \rightarrow V$ denote the linear transformation such that for $0 \leq i \leq d$, $U_i$ is an eigenspace of $K$ with eigenvalue $q^{2i-d}$.
\end{definition}

\begin{remark}
\rm
Combining (\ref{AA^*U1}) and Definition \ref{def:K} we have
\begin{eqnarray}
\label{AKU}
(A-K)U_i \subseteq U_{i+1} \qquad \qquad (0 \leq i \leq d), \\
\label{A^*KK^-1U}
(A^*-K-cK^{-1})U_i \subseteq U_{i-1} \qquad \qquad (0 \leq i \leq d).
\end{eqnarray}
\end{remark}

\noindent
The goal for the remainder of this section is to prove a number of relations between the linear transformations $A,A^*,B,K,K^{-1}$ which will be used in Section 13.

\begin{lemma}
\label{thm:AA^*K}
With reference to Assumption \ref{assum} and Definition \ref{def:K} we have
\begin{eqnarray}
\label{AA^*K1}
\frac{qK^{-1}A-q^{-1}AK^{-1}}{q-q^{-1}} &=& I, \\
\label{AA^*K2}
\frac{qKA^*-q^{-1}A^*K}{q-q^{-1}} &=& K^2 + c \, I.
\end{eqnarray}
\end{lemma}

\noindent
{\it Proof:} Recall that  $\lbrace U_i \rbrace _{i=0}^d$ is a decomposition of $V$.  Combining Definition \ref{def:K}, (\ref{AKU}), and Lemma \ref{thm:bidiag1} we obtain (\ref{AA^*K1}).  Combining Definition \ref{def:K}, (\ref{A^*KK^-1U}), and Lemma \ref{thm:bidiag3} we obtain (\ref{AA^*K2}).
\hfill $\Box$ \\

\begin{lemma}
With reference to Lemma \ref{def:U_i}, Definition \ref{def:B}, and Definition \ref{def:K} we have
\begin{eqnarray}
\label{A^*BU1}
(B-K)U_i \subseteq U_{i-1} \qquad \qquad (0 \leq i \leq d), \\
\label{A^*BU2}
(A^*-B-cK^{-1})U_i \subseteq U_{i-1} \qquad \qquad (0 \leq i \leq d).
\end{eqnarray}
\end{lemma}

\noindent
{\it Proof:} First we show (\ref{A^*BU1}).  Using Lemma \ref{thm:bidiag1} and (\ref{AA^*B1}) we have
\begin{eqnarray}
\label{BKrel2}
(B-q^{d-2i}I)V_i \subseteq V_{i-1} \qquad \qquad (0 \leq i \leq d).
\end{eqnarray}
We have
\[ \begin{array}{cccc}
(B-K) U_i &=&
(B-q^{2i-d}I) U_i  \hfill&\hskip 1 truein  \quad (\mbox{by Definition \ref{def:K}}) \\
&\subseteq& (B-q^{2i-d}I) (U_0+ \cdots + U_i) \hfill&
\\
&=& (B-q^{2i-d}I) (W_0+ \cdots + W_i) \hfill &\hskip 1 truein  \quad (\mbox{by (\ref{AA^*U2}), (\ref{AA^*W2})})
\\
&\subseteq& W_0+ \cdots + W_{i-1} \hfill&\hskip 1 truein  \quad (\mbox{by Definition \ref{def:B}})
\\
&=& U_0+ \cdots + U_{i-1} \hfill&\hskip 1 truein \quad (\mbox{by (\ref{AA^*U2}), (\ref{AA^*W2})}) \end{array}\] \\
and also \\
\[ \begin{array}{cccc}
(B-K) U_i &=&
(B-q^{2i-d}I) U_i  \hfill&\hskip 1 truein  \quad (\mbox{by Definition \ref{def:K}}) \\
&\subseteq& (B-q^{2i-d}I) (U_i+ \cdots + U_d)\hfill&
\\
&=&(B-q^{2i-d}I) (V_i+ \cdots + V_d)\hfill&\hskip 1 truein \quad (\mbox{by (\ref{AA^*U2})})
\\
&\subseteq& V_{i-1}+ \cdots + V_d \hfill&\hskip 1 truein\quad (\mbox{by (\ref{BKrel2})})
\\
&=&
U_{i-1}+ \cdots + U_d\hfill&\hskip 1 truein\quad (\mbox{by (\ref{AA^*U2})}).
\end{array}\]
Using this and since $\lbrace U_i \rbrace _{i=0}^d$ is a decomposition of $V$ we have (\ref{A^*BU1}). Combining (\ref{A^*KK^-1U}) and (\ref{A^*BU1}) we obtain (\ref{A^*BU2}).
\hfill $\Box$ \\

\begin{lemma}
\label{thm:BKrel}
With reference to Definition \ref{def:B} and Definition \ref{def:K} we have
\begin{eqnarray}
\label{BKrel1}
\frac{qBK^{-1}-q^{-1}K^{-1}B}{q-q^{-1}} = I.
\end{eqnarray}
\end{lemma}

\noindent
{\it Proof:} Recall that  $\lbrace U_i \rbrace _{i=0}^d$ is a decomposition of $V$.  Combining Lemma \ref{thm:bidiag2} and (\ref{A^*BU1}) we obtain (\ref{BKrel1}).
\hfill $\Box$ \\

\begin{lemma}
\label{thm:L'LR}
With reference to Assumption \ref{assum}, Definition \ref{def:B} and Definition \ref{def:K} for an integer $j \geq 1$ we have
\begin{eqnarray}
\label{L'LR1}
\frac{q^{j} (B-K)^{j} (A^*-K-cK^{-1}) - q^{-j} (A^*-K-cK^{-1}) (B-K)^{j}}{q^{j}-q^{-j}} = (B-K)^{j+1},
\end{eqnarray}
\begin{eqnarray}
\label{L'LR2}
\frac{q^{j} (A-K)^{j} (B-K) - q^{-j} (B-K) (A-K)^{j}}{q^{j}-q^{-j}} = -(q^{2-2j}K^2-I) (A-K)^{j-1}.
\end{eqnarray}
\end{lemma}

\noindent
{\it Proof:} First we show (\ref{L'LR1}) by induction on $j$.  Multiplying out the left hand side of (\ref{L'LR1}) with $j=1$ and simplifying the result using (\ref{AA^*B2}), (\ref{AA^*K2}), and (\ref{BKrel1}) we obtain the right hand side of (\ref{L'LR1}) with $j=1$.  To prove (\ref{L'LR1}) for $j \geq 2$ note that (\ref{L'LR1}) is equivalent to
\begin{eqnarray*}
(B-K)^{j} (A^*-K-cK^{-1}) =  q^{-2j} (A^*-K-cK^{-1}) (B-K)^{j} + q^{-j} (q^{j}-q^{-j}) (B-K)^{j+1}.
\end{eqnarray*}
This is shown by a routine induction argument using the $j=1$ case.  We now show (\ref{L'LR2}) by induction on $j$.  Multiplying out the left hand side of (\ref{L'LR2}) with $j=1$ and simplifying the result using (\ref{AA^*B1}), (\ref{AA^*K1}), and (\ref{BKrel1}) we obtain the right hand side of (\ref{L'LR2}) with $j=1$.  Note that (\ref{AA^*K1}) is equivalent to
\begin{eqnarray}
\label{L'LR3}
(A-K)K=q^{-2}K(A-K).
\end{eqnarray}
To prove (\ref{L'LR2}) for $j \geq 2$ note that (\ref{L'LR2}) is equivalent to  \begin{eqnarray*}
(B-K)(A-K)^{j} =  q^{2j} (A-K)^{j} (B-K) + q^{2-j} (q^{j}-q^{-j}) (K^2-q^{2j-2}I) (A-K)^{j-1}.
\end{eqnarray*}
This is shown by a routine induction argument using the $j=1$ case and (\ref{L'LR3}).
\hfill $\Box$ \\

\section{The projections $E_i, \widetilde{E}^*_i$}

\begin{definition}
\label{def:EiE^*i}
\rm
With reference to Assumption \ref{assum} and Definition \ref{def:V^*_itilde} we define the following linear transformations.
\begin{enumerate}
\item For $0 \leq i \leq d$, we let $E_i$ denote the linear transformation on $V$ satisfying both
\begin{center}
$(E_i - I)V_i=0$, \\
$E_iV_j=0 \,\,\,\,  \mbox{if} \,\,\,\, j \neq i, \qquad (0 \leq j \leq d)$.
\end{center}
\item For $0 \leq i \leq d$, we let $\widetilde{E}^*_i$ denote the linear transformation on $V$ satisfying both
\begin{center}
$(\widetilde{E}^*_i - I)\widetilde{V}^*_i=0$, \\
$\widetilde{E}^*_i\widetilde{V}^*_j=0 \,\,\,\,  \mbox{if} \,\,\,\, j \neq i, \qquad (0 \leq j \leq d)$.
\end{center}
\end{enumerate}
In other words, $E_i$ (resp. $\widetilde{E}^*_i$) is the projection map from $V$ onto $V_i$ (resp. $\widetilde{V}^*_i$).
\end{definition}

\begin{lemma}
\label{thm:EE^*AA^*}
With reference to Assumption \ref{assum}, Definition \ref{def:A^*tilde}, and Definition \ref{def:EiE^*i}, for $0 \leq i \leq d$ we have
\begin{eqnarray}
E_i &=& \prod_{{0 \leq j \leq d}\atop {j\not=i}}
\frac{A-q^{2j-d}I}{q^{2i-d}-q^{2j-d}},
\label{EE^*AA^*1}
\\
\widetilde{E}^*_i &=& \prod_{{0 \leq j \leq d}\atop {j\not=i}}
\frac{\widetilde{A}^*-q^{d-2j}I}{q^{d-2i}-q^{d-2j}}.
\label{EE^*AA^*2}
\end{eqnarray}
\end{lemma}

\noindent {\it Proof:}
Concerning (\ref{EE^*AA^*1}), let $E'_i$ denote the expression on the right in that line.  Using Assumption \ref{assum} we find $(E'_i-I)V_i=0$ and $E'_iV_j=0$ $(0 \leq j \leq d, \; j\not=i)$.  By this and Definition \ref{def:EiE^*i}(i) we find $E_i=E'_i$.  We have now proved (\ref{EE^*AA^*1}).  The proof of (\ref{EE^*AA^*2}) is similar.
\hfill $\Box $ \\

\begin{lemma}
\label{thm:bij}
With reference to Assumption \ref{assum}, Remark \ref{G_iformulas}, and Definition \ref{def:EiE^*i}(i) the following holds for $0 \leq i \leq d$:
The linear transformations
\begin{eqnarray*}
{{W_{d-i}\quad \rightarrow \quad  V_i}\atop {w \quad \rightarrow \quad E_iw}}
\qquad \qquad \qquad
{{V_i\quad \rightarrow \quad W_{d-i}}\atop {v \quad \; \rightarrow \quad G_{d-i}v}}
\end{eqnarray*}
are bijections, and  moreover, they are inverses.
\end{lemma}

\noindent
{\it Proof:}  It suffices to show $G_{d-i}E_i-I$ vanishes on $W_{d-i}$ and $E_iG_{d-i}-I$ vanishes on $V_i$.  We will use the following notation.  Recall by (\ref{AA^*W2}) that for $0 \leq j \leq d$, $W_{d-j}+ \cdots + W_d = V_0+ \cdots + V_j$;  let $Z_j$ denote this common sum.  We set $Z_{-1}=0$.  By the construction $Z_i=W_{d-i} + Z_{i-1}$ (direct sum) and $Z_i=V_i + Z_{i-1}$ (direct sum).  Also $(I-G_{d-i})Z_i=Z_{i-1}$ and $(I-E_i)Z_i=Z_{i-1}$.  We now show $G_{d-i}E_i-I$ vanishes on $W_{d-i}$.  Pick $w \in W_{d-i}$.  Using $G_{d-i}E_i-I=(G_{d-i}-I)E_i + E_i-I$ and our preliminary comments we routinely find $(G_{d-i}E_i-I)w \in Z_{i-1}$.  But $(G_{d-i}E_i-I)w \in W_{d-i}$ by construction and $W_{d-i}\cap Z_{i-1}=0$ so $(G_{d-i}E_i-I)w=0$.  We now show $E_iG_{d-i}-I$ vanishes on $V_i$.  Pick $v \in V_i$.  Using $E_iG_{d-i}-I=(E_i-I)G_{d-i}+G_{d-i}-I$ and our preliminary comments we routinely find $(E_iG_{d-i}-I)v \in Z_{i-1}$.  But $(E_iG_{d-i}-I)v \in V_i$ by construction and $V_i\cap Z_{i-1}=0$ so $(E_iG_{d-i}-I)v=0$.  We have now shown $G_{d-i}E_i-I$ vanishes on $W_{d-i}$ and $E_iG_{d-i}-I$ vanishes on $V_i$.  Consequently the given maps are inverses.  Each of these maps has an inverse and is therefore a bijection.
\hfill $\Box $ \\

\section{How $\widetilde{E}_0^*, E_d, P$ are related}

\noindent
The goal of this section is to prove the following theorem which will be used in the proof of Theorem \ref{thm:PVnot0}.

\begin{theorem}
\label{thm:PVformula}
With reference to Assumption \ref{assum}, Lemma \ref{def:U_i}, Definition \ref{def:PV} and Definition \ref{def:EiE^*i}, for $u \in U_0$ we have
\begin{eqnarray}
\widetilde{E}^*_0E_du = c^{-d} \, q^{2d(1-d)} \, P( q^{2d-2} (q-q^{-1})^{-2} ) \, u.
\label{mainformula}
\end{eqnarray}
\end{theorem}

\noindent
Before we prove Theorem \ref{thm:PVformula} we develop some notation and prove some preliminary lemmas.  \\

\noindent
With reference to (\ref{def:brackfacn}) for integers $n,m$ with $n \geq 0$ and $0 \leq m \leq n$ we define
\begin{eqnarray}
\label{bracket}
\left[ n \atop m \right]  = \frac{\lbrack n \rbrack !}{\lbrack m \rbrack ! \, \lbrack n-m \rbrack !}.
\end{eqnarray}

\begin{lemma}
\rm
For integers $n,m$ with $n \geq 1$ and $0 \leq m \leq n-1$ we have
\begin{eqnarray}
\label{brakformula1}
\left[ n-1 \atop m \right] \, + \, q^{n} \left[ n-1 \atop m-1 \right] = q^{m} \left[ n \atop m \right], \\
\label{brakformula2}
\left[ n-1 \atop m \right] \, + \, q^{-n} \left[ n-1 \atop m-1 \right] = q^{-m} \left[ n \atop m \right].
\end{eqnarray}
\end{lemma}

\noindent
{\it Proof:}  Immediate from (\ref{def:brackn}), (\ref{def:brackfacn}) and (\ref{bracket}).
\hfill $\Box$ \\

\noindent
The following two lemmas provide key formulas to be used in the proof of Theorem \ref{thm:PVformula}.

\begin{lemma}
\label{thm:BpassA^*}
With reference to Assumption \ref{assum}, Definition \ref{def:B} and Definition \ref{def:K} for an integer $i \geq 0$ we have
\begin{eqnarray}
(A^*-B-cK^{-1})^{i} \, = \, \sum_{j=0}^{i} (-1)^{j} \, q^{j-ji} \left[ i \atop j \right] (A^*-K-cK^{-1})^{i-j} \, (B-K)^{j}.
\label{BpassA^*1}
\end{eqnarray}
\end{lemma}

\noindent
{\it Proof:}  We prove (\ref{BpassA^*1}) by induction on $i$. For $i=0$ both sides of (\ref{BpassA^*1}) equal $I$.  Now let $i \geq 1$.  Abbreviate $\Delta=A^*-K-cK^{-1}$ and $\Gamma=B-K$.  We have
\[ \begin{array}{cclr}
(A^*-B-cK^{-1})^{i} &=& \medskip
(A^*-B-cK^{-1})^{i-1} (\Delta - \Gamma) \hfill&
\\ \medskip
&=& \sum_{j=0}^{i-1} (-1)^{j} q^{2j-ji} \left[ i-1 \atop j \right] \Delta^{i-j-1} (\Gamma^j \Delta - \Gamma^{j+1}) \hfill&  \,\,\,\,\,\,\,\,\,\,\,\,\,\,\,\,\,\,\,\,\,\, (\mbox{by induction})
\\ \medskip
&=& \sum_{j=0}^{i-1} (-1)^{j} q^{-ji} \left[ i-1 \atop j \right] \Delta^{i-j-1} (\Delta \Gamma^j - \Gamma^{j+1}) \hfill&   (\mbox{by (\ref{L'LR1})})
\\ \medskip
&=& \sum_{j=0}^{i-1} (-1)^{j} q^{-ji} \left[ i-1 \atop j \right] \Delta^{i-j} \Gamma^j \, \\ && \,\,\,\,\,\,\,\,\,\,\,\, + \,\, \sum_{j=1}^{i} (-1)^{j} q^{i-ji} \left[ i-1 \atop j-1 \right] \Delta^{i-j} \Gamma^j \hfill&
\\ \medskip
&=& \sum_{j=0}^{i} (-1)^{j} q^{j-ji} \left[ i \atop j \right] \Delta^{i-j} \Gamma^j \hfill&   (\mbox{by (\ref{brakformula1})}).
\end{array}\]
\hfill $\Box$ \\

\begin{lemma}
\label{thm:BpassA}
Fix an integer $i \geq 1$.  With reference to Definition \ref{def:K}, for integers $\mu, \nu \geq 0$, define a polynomial $f_{\mu, \nu} \in \K[K^2]$ by                 $f_{\mu, \nu}=\prod_{s=0}^{\mu -1}(K^2-q^{2i-2s-2 \nu}I)$.  With reference to Assumption \ref{assum} and Definition \ref{def:B} for $1 \leq j \leq i$ we have
\begin{eqnarray}
(B-K)^{j} \, (A-K)^{i} = \sum_{h=0}^{j} \, q^{M_{h,i,j}} \, N_{h,i} \, \left[ j \atop h \right]  \, f_{h,j} \,\, (A-K)^{i-h} \, (B-K)^{j-h}
\label{BpassA1}
\end{eqnarray}
where $M_{h,i,j}={(h/2)(3h-1) + hj - 3hi + 2ij}$ and $N_{h,i} = \left[ i \atop h \right] \lbrack h \rbrack ! \, (q-q^{-1})^h$.
\end{lemma}

\noindent
{\it Proof:} Before we prove (\ref{BpassA1}) we have a comment.  Observe that (\ref{BKrel1}) is equivalent to
\begin{eqnarray}
(B - K) K = q^{2} K (B - K).
\label{BpassA2}
\end{eqnarray}
To prove (\ref{BpassA1}) we let $i$ be given and use induction  on $j$.  For $j=1$ (\ref{BpassA1}) is equivalent to (\ref{L'LR2}).  Now let $j \geq 2$.  Abbreviate $\Delta=A-K$ and $\Gamma=B-K$.  We have
\[ \begin{array}{cccc}
\Gamma^j \Delta^i &=& \medskip
\Gamma \Gamma^{j-1} \Delta^i \hfill&
\\ \medskip
&=& \Gamma \, \sum_{h=0}^{j-1} \, q^{M_{h,i,j-1}} \, N_{h,i} \, \left[ j-1 \atop h \right] \, \, f_{h,j-1} \, \Delta^{i-h} \, \Gamma^{j-h-1} \,\,\,\,\,\,\,\,\,\,\,\,\,\,\,\,\,\,\,\,\,\,\,\,\,\,\,\,\,\,\, {\rm (by \,\, induction)} \hfill&
\\ \medskip
&=& \sum_{h=0}^{j-1} \, q^{M_{h,i,j}} \, q^{-h-2i} \, N_{h,i} \, \left[ j-1 \atop h \right] \, q^{4h} \, f_{h,j+1} \, \Gamma \, \Delta^{i-h} \, \Gamma^{j-h-1} \,\,\,\,\,\,\,\,\,\,\,\,\, {\rm (by \,\, (\ref{BpassA2}))} \hfill&
\\ \medskip
&=& \sum_{h=0}^{j-1} \, q^{M_{h,i,j}} \, q^{3h-2i} \, N_{h,i} \, \left[ j-1 \atop h \right]  \, f_{h,j+1} \, q^{2i-2h} \, \Delta^{i-h} \, \Gamma^{j-h} \,\,\,\,\,\,\,\,\,\,\,\,\,\,\,\,\,\, {\rm (by \,\, (\ref{L'LR2}))} \hfill&
\\ \medskip
& & \,\,\,\,\,\,\,\, + \, \sum_{h=0}^{j-1} \, q^{M_{h+1,i,j}} \, q^{h+1-j} \, N_{h+1,i} \, \left[ j-1 \atop h \right] \, f_{h,j+1} \, (K^2-q^{2i-2h-2}I) \, \Delta^{i-h-1} \, \Gamma^{j-h-1}  \hfill&
\\ \medskip
&=& \sum_{h=0}^{j-1} \, q^{M_{h,i,j}} \, q^{h} \, N_{h,i} \, \left[ j-1 \atop h \right]  \, f_{h-1,j+1} \, (K^2 - q^{2i-2j-2h}I) \, \Delta^{i-h} \, \Gamma^{j-h} \hfill&
\\ \medskip
& & \,\,\,\,\,\,\,\, + \, \sum_{h=1}^{j} \, q^{M_{h,i,j}} \, q^{h-j} \, N_{h,i} \, \left[ j-1 \atop h-1 \right] \, f_{h-1,j+1} \, (K^2-q^{2i-2h}I) \, \Delta^{i-h} \, \Gamma^{j-h} \hfill&
\\ \medskip
&=& \sum_{h=0}^{j} \, q^{M_{h,i,j}} \, N_{h,i} \, \left[ j \atop h \right] \, f_{h-1,j+1} \, (K^2 - q^{2i-2j}I) \, \Delta^{i-h} \, \Gamma^{j-h} \,\,\,\,\,\,\,\,\,\,\, {\rm (by \,\, (\ref{brakformula1}), (\ref{brakformula2}))} \hfill&
\\ \medskip
&=& \sum_{h=0}^{j} \, q^{M_{h,i,j}} \, N_{h,i} \, \left[ j \atop h \right] \, f_{h,j} \, \Delta^{i-h} \, \Gamma^{j-h}. \hfill&
\end{array}\]
\hfill $\Box$ \\

\noindent
We are now ready to prove Theorem \ref{thm:PVformula}. \\

\noindent
{\it Proof of Theorem \ref{thm:PVformula}:}  Let $u \in U_0$.  Using Definition \ref{def:A^*tilde} and Lemma \ref{thm:EE^*AA^*} we have
\begin{eqnarray}
\widetilde{E}^*_0E_du = c^{-d} \, q^{d-d^2} \, (q-q^{-1})^{-2d} \, \lbrack d \rbrack !^{-2} \, \prod_{j=1}^{d} (A^*-B-cq^{d-2j}I) \prod_{j=0}^{d-1} (A-q^{2j-d}I) u.
\label{PVformula1}
\end{eqnarray}
Applying Definition \ref{def:K}, (\ref{AKU}), and (\ref{A^*BU2}) to (\ref{PVformula1}) we have
\begin{eqnarray}
\widetilde{E}^*_0E_du = c^{-d} \, q^{d-d^2} \, (q-q^{-1})^{-2d} \, \lbrack d \rbrack !^{-2} \, (A^*-B-cK^{-1})^d \, (A-K)^d u.
\label{PVformula2}
\end{eqnarray}
We now express the right hand side of (\ref{PVformula2}) in terms of the maps $R,L$ from Definition \ref{def:RL}.  By (\ref{A^*BU1}) we have $(B-K)u=0$.  Using Definition \ref{def:Fi}, (\ref{RL1}), and Definition \ref{def:K} we have $R=A-K$ and $L=A^*-K-cK^{-1}$.  Using Lemma \ref{thm:BpassA^*}, Lemma \ref{thm:BpassA}, and the previous two sentences we have
\begin{eqnarray}
(A^*-B-cK^{-1})^d \, (A-K)^d u = \sum_{j=0}^{d} C_{j} \,\, L^{d-j} \, \prod_{s=0}^{j -1}(K^2-q^{2d-2s-2j}I) \, R^{d-j} u
\label{PVformula3}
\end{eqnarray}
where $C_{j}=(-1)^{j} \, q^{(j/2)(5j+1)-2dj} \, \left[ d \atop j \right] ^2 \, \lbrack j \rbrack ! \, (q-q^{-1})^{j}$.  \\ Using Definition \ref{def:RL} we have $R^{d-j}u \in U_{d-j}$ for $0 \leq j \leq d$.  So by Definition \ref{def:K} \\ $(K^2-q^{2d-4j}I)R^{d-j}u=0$.  Using this on the right hand side of (\ref{PVformula3}) and simplifying the result we have
\begin{eqnarray}
(A^*-B-cK^{-1})^d \, (A-K)^d u = \sum_{j=0}^{d} q^{j-j^2} \, \lbrack d \rbrack ! ^2 \, \lbrack d-j \rbrack !^{-2} \, (q-q^{-1})^{2j} \, L^{d-j} \, R^{d-j} u.
\label{PVformula4}
\end{eqnarray}
Changing the index of summation in (\ref{PVformula4}) by letting $j=d-t$, using Definition \ref{def:zeta}, and simplifying the result we have
\begin{eqnarray*}
(A^*-B-cK^{-1})^d (A-K)^d u = q^{d-d^2} \, \lbrack d \rbrack ! ^2 \, (q-q^{-1})^{2d} \sum_{t=0}^{d} q^{t(1-t)} \, \lbrack t \rbrack !^{-2} \, (q^{2d-2} \, (q-q^{-1})^{-2})^t \, \zeta_{t} u.
\label{PVformula5}
\end{eqnarray*}
Combining the previous line with (\ref{PVformula2}) and using Definition \ref{def:PV} we obtain (\ref{mainformula}).
\hfill $\Box$ \\

\section{The raising/lowering maps revisited}

\noindent
In this section we prove a number of relations between $r,l$ from Definition \ref{def:W_iandfriends} and $B$ from Definition \ref{def:B}.  These relations will help to motivate the next section.

\begin{lemma}
\label{thm:AA^*RLrl}
With reference to Assumption \ref{assum}, Definition \ref{def:W_iandfriends}, and Definition \ref{def:B} we have
\begin{enumerate}
\item $r=A-B^{-1}$,
\item $l=A^*-B-cB^{-1}$.
\end{enumerate}
\end{lemma}

\noindent
{\it Proof:} (i) Recall by Assumption \ref{assum} that $\theta_i=q^{2i-d}$ for $0 \leq i \leq d$.  Using Remark \ref{G_iformulas} and Definition \ref{def:B} we find $\sum_{i=0}^{d} \theta_{d-i} G_i =B^{-1}$  Using this and (\ref{rl}) we obtain the desired result. \\
(ii) Similar to (i).
\hfill $\Box$ \\

\begin{lemma}
\label{thm:RLrlBK}
With reference to Definition \ref{def:W_iandfriends} and Definition \ref{def:B} we have
\begin{enumerate}
\item $Br=q^2rB$,
\item $Bl=q^{-2}lB$.
\end{enumerate}
\end{lemma}

\noindent
{\it Proof:} (i) Recall by Definition \ref{def:W_iandfriends} that  $\lbrace W_i \rbrace _{i=0}^d$ is a decomposition of $V$.  So it suffices to show $Br-q^2rB$ vanishes on $W_i$ for $0 \leq i \leq d$. Let $i$ be given and let $w \in W_i$.  Using Remark \ref{rlformulas} and Definition \ref{def:B} we find $rw$ is an eigenvector for $B$ with eigenvalue $q^{2i+2-d}$.  From this we find $(Br-q^2rB)w=0$ and the result follows.  \\
(ii) Similar to (i).
\hfill $\Box$ \\

\begin{lemma}
\label{thm:RLrlrel}
With reference to Definition \ref{def:W_iandfriends} and Definition \ref{def:B} we have
\begin{enumerate}
\item $r^3l - \lbrack 3 \rbrack r^2lr + \lbrack 3 \rbrack rlr^2 - lr^3=q^{-4} \, (q-q^{-1})^3 \, \lbrack 3 \rbrack ! \, r^2B^{-2}$,
\item $rl^3 - \lbrack 3 \rbrack lrl^2 + \lbrack 3 \rbrack l^2rl - l^3r = q^{-4} \, (q-q^{-1})^3 \, \lbrack 3 \rbrack ! \, B^{-2}l^2$.
\end{enumerate}
\end{lemma}

\noindent
{\it Proof:}
By Lemma \ref{thm:AA^*RLrl} we have $A=r+B^{-1}$ and $A^*= l + B + cB^{-1}$.  Substituting these into Lemma \ref{thm:AA^*rel}(i)(ii) and simplifying using Lemma \ref{thm:RLrlBK} we obtain the desired result.
\hfill $\Box$ \\

\section{The algebra $A_q {(\alpha)}$}

\noindent
Motivated by Lemma \ref{thm:RLrlBK} and Lemma \ref{thm:RLrlrel} we define an algebra $A_q(\alpha)$.  We find a spanning set for $A_q(\alpha)$ that will be used in the proof of Lemma \ref{thm:t=d} in the next section.

\begin{definition}
\rm
\label{def:A'qalgebra}
Fix a scalar $\alpha \in \K$.  Let $A_q(\alpha)$ denote the unital associative $\K$-algebra defined by generators $x,y,z,z^{-1}$ subject to the relations
\begin{eqnarray}
\label{A'q1}
zz^{-1}&=&1\,\,=\,\,z^{-1}z, \\
\label{A'q2}
zx&=&q^2xz, \\
\label{A'q3}
zy&=&q^{-2}yz,
\end{eqnarray}
\vspace{-.4 in}
\begin{eqnarray}
\label{A'q4}
x^3y - \lbrack 3 \rbrack x^2yx + \lbrack 3 \rbrack xyx^2 - yx^3&=& \alpha \, x^2z^{-2}, \\
\label{A'q5}
xy^3 - \lbrack 3 \rbrack yxy^2 + \lbrack 3 \rbrack y^2xy  - y^3x&=&\alpha \, z^{-2}y^{2}.
\end{eqnarray}
\end{definition}

\begin{note}
\rm
In the case $\alpha=0$ we note that the algebra $A_q(0)$ is the algebra $\mathcal{B}$ from \cite[Definition~1.10]{Benkart04}.  The algebra $A_q(\alpha)$ is a special case of a more general algebra currently being studied called the augmented tridiagonal algebra.  The proof of Theorem \ref{thm:irredbasisforAq} involves an argument which is adapted from [T. Ito and P. Terwilliger, {\it The augmented tridiagonal algebra}, in preparation].
\end{note}

\noindent
Before we display a spanning set for $A_q(\alpha)$ we have a number of preliminary comments.

\medskip

\noindent
For the moment we view $x,y$ as formal symbols and let $F$ denote the free unital associative $\K$-algebra on $x,y$.

\begin{definition}
\rm
\label{def:redword}
By a {\it word in $F$} we mean an element of $F$ of the form $a_1 a_2 \cdots a_n$ where $n$ is a nonnegative integer and $a_i \in \lbrace x,y \rbrace$ for $1 \leq i \leq n$.  We call $n$ the {\it length} of $a_1 a_2 \cdots a_n$.  We interpret the word of length 0 as the identity element of $F$.  We say this word is {\it trivial}.  Observe $F = \sum_{n=0}^{\infty} F_n$ (direct sum) where $F_n$ denotes the subspace of $F$ spanned by all the words of length $n$.  Moreover, $F_n F_m = F_{n+m}$.
\end{definition}

\begin{definition}
\rm
Let $a_1 a_2 \cdots a_n$ denote a word in $F$.  Observe there exists a unique sequence $(i_1, i_2, \ldots, i_s)$ of positive integers such that $a_1 a_2 \cdots a_n$ is one of $x^{i_1}y^{i_2}x^{i_3} \cdots y^{i_s}$ or $x^{i_1}y^{i_2}x^{i_3} \cdots x^{i_s}$ or $y^{i_1}x^{i_2}y^{i_3} \cdots x^{i_s}$ or $y^{i_1}x^{i_2}y^{i_3} \cdots y^{i_s}$.  We call the sequence $(i_1, i_2, \ldots, i_s)$ the {\it signature} of $a_1a_2 \cdots a_n$.
\end{definition}

\begin{example}
\rm
Each of the words $yx^2y^2x$, $xy^2x^2y$ has signature $(1,2,2,1)$.
\end{example}

\begin{definition}
\label{def:irredword}
\rm
Let $a_1 a_2 \cdots a_n$ denote a word in $F$ and let $(i_1, i_2, \ldots, i_s)$ denote the corresponding signature.  We say $a_1 a_2 \cdots a_n$ is {\it reducible} whenever there exists an integer $\eta$ $(2 \leq \eta \leq s-1)$ such that $i_{\eta-1} \geq i_\eta < i_{\eta+1}$.  We say a word in $F$ is {\it irreducible} whenever it is not reducible.
\end{definition}

\begin{example}
\rm
A word in $F$ of length less than 4 is irreducible.  The only reducible words in $F$ of length 4 are $xyx^2$ and $yxy^2$.
\end{example}

\noindent
In the following lemma we give a necessary and sufficient condition for a given nontrivial word in $F$ to be irreducible.

\begin{lemma}
\label{thm:irredform}
Let $a_1 a_2 \cdots a_n$ denote a nontrivial word in $F$ and let $(i_1, i_2, \ldots, i_s)$ denote the corresponding signature.  Then the following are equivalent:
\begin{enumerate}
\item[{\rm (i)}]
The word $a_1 a_2 \cdots a_n$ is irreducible.
\item[{\rm (ii)}]
There exists an integer $t$ $(1 \leq t \leq s)$ such that
\begin{eqnarray*}
i_1 < i_2 < \cdots < i_t \geq i_{t+1} \geq i_{t+2} \geq \cdots \geq i_{s-1} \geq i_s.
\end{eqnarray*}
\end{enumerate}
\end{lemma}

\noindent
{\it Proof:}  Immediate from Definition \ref{def:irredword}.
\hfill $\Box$ \\

\noindent
For the moment we view $x,y,z,z^{-1}$ as formal symbols and let ${\mathcal F}$ denote the free unital associative $\K$-algebra on $x,y,z,z^{-1}$.  We identify $F$ with the subalgebra of ${\mathcal F}$ generated by $x,y$.

\medskip

\noindent
We now view $A_q(\alpha)$ as a vector space over $\K$ and display a spanning set.

\begin{theorem}
\label{thm:irredbasisforAq}
 Let $\pi : {\mathcal F} \rightarrow A_q(\alpha)$ denote the canonical quotient map.  Consider the following elements in ${\mathcal F}$:
\begin{center}
$wz^j$, \quad $w$ is an irreducible word in $F$, \quad $j \in \Z$.
\end{center}
Then $A_q(\alpha)$ is spanned by the images of the above elements under $\pi$.
\end{theorem}

\noindent
To prove Theorem \ref{thm:irredbasisforAq} we will need the following two lemmas and definition.

\begin{lemma}
\label{thm:shaperesult}
Let $\Omega$ denote the subspace of $F$ spanned by all the irreducible words.  Let $\Lambda$ denote the two sided ideal of $F$ generated by
\begin{eqnarray}
\label{gen1}
x^3y - \lbrack 3 \rbrack x^2yx + \lbrack 3 \rbrack xyx^2 - yx^3, \\
\label{gen2}
xy^3 - \lbrack 3 \rbrack yxy^2 + \lbrack 3 \rbrack y^2xy  - y^3x.
\end{eqnarray}
For an integer $n \geq 0$ let $\Omega_n = \Omega \cap F_n$ and $\Lambda_n = \Lambda \cap F_n$.  Then the following (i)--(iv) hold:
\begin{enumerate}
\item $F = \Omega + \Lambda$  \qquad (direct sum),
\item $\Omega = \sum_{n=0}^{\infty} \Omega_n$ \qquad (direct sum),
\item $\Lambda = \sum_{n=0}^{\infty} \Lambda_n$ \qquad (direct sum),
\item $F_n = \Omega_n + \Lambda_n$ \qquad (direct sum) \qquad $0 \leq n < \infty $.
\end{enumerate}
\end{lemma}

\noindent
{\it Proof:} (i) View the $\K$-algebra $F/ \Lambda$ as a vector space over $\K$.  By \cite[Theorem 2.29]{ItoTer04} $F/ \Lambda$ has a basis consisting of the images of the irreducible words in $F$ under the canonical quotient map $F \rightarrow F/ \Lambda$.  The result follows immediately from this.  \\
(ii) The words in $F$ form a basis for $F$. \\
(iii) The generators (\ref{gen1}), (\ref{gen2}) of $\Lambda$ are in $F_4$.  \\
(iv) Recall $F = \sum_{n=0}^{\infty} F_n$ (direct sum).  Combining this with (i)--(iii) above we obtain the desired result.
\hfill $\Box$ \\

\begin{lemma}
\label{thm:pilambda}
We have $\Lambda_n = 0$ for $n \leq 3$.  Also
\begin{eqnarray}
\label{pilam}
\pi (\Lambda_n) \, \subseteq \, \pi (F_{n-2} \, z^{-2}) \qquad \qquad n \geq 4
\end{eqnarray}
where $\pi : {\mathcal F} \rightarrow A_q(\alpha)$ is the canonical quotient map.
\end{lemma}

\noindent
{\it Proof:}  The first assertion follows since the generators (\ref{gen1}), (\ref{gen2}) of $\Lambda$ are in $F_4$.  For $n \geq 4$ we have by construction that
\begin{eqnarray*}
\Lambda_n \,\, = \,\, \sum_{i,j} F_i \, (x^3y - \lbrack 3 \rbrack x^2yx + \lbrack 3 \rbrack xyx^2 - yx^3) \, F_j
\end{eqnarray*}
\vspace{-.1in}
\begin{eqnarray*}
\qquad \qquad \qquad + \sum_{i,j} F_i \, (xy^3 - \lbrack 3 \rbrack yxy^2 + \lbrack 3 \rbrack y^2xy  - y^3x) \, F_j
\end{eqnarray*}
where each sum is over all nonnegative integers $i,j$ such that $i+j = n-4$.  Applying $\pi$ and using (\ref{A'q2})--(\ref{A'q5}) we have
\begin{eqnarray}
\label{lambformula}
\pi (\Lambda_n) \,\, = \,\, \sum_{i,j} \pi ( F_i \, x^2 \, F_j \, z^{-2}) \,\, + \,\, \sum_{i,j} \pi ( F_i \, y^2 \, F_j \, z^{-2})
\end{eqnarray}
where each sum is over all nonnegative integers $i,j$ such that $i+j = n-4$.  For all such $i,j$ we have $F_i \, x^2 \, F_j \subseteq F_{n-2}$ and $F_i \, y^2 \, F_j \subseteq F_{n-2}$.  Simplifying (\ref{lambformula}) using this we obtain (\ref{pilam}).
\hfill $\Box$ \\

\begin{definition}
\rm
\label{def:wordsinscriptF}
By a {\it word in ${\mathcal F}$} we mean an element of ${\mathcal F}$ of the form $a_1 a_2 \cdots a_n$ where $n$ is a nonnegative integer and $a_i \in \lbrace x,y,z,z^{-1} \rbrace$ for $1 \leq i \leq n$.  By the {\it $(x,y)$-length} of $a_1 a_2 \cdots a_n$ we mean the number of $x$'s plus the number of $y$'s in $a_1 a_2 \cdots a_n$.
\end{definition}

\noindent
We are now ready to prove Theorem \ref{thm:irredbasisforAq}. \\

\noindent
{\it Proof of Theorem \ref{thm:irredbasisforAq}:}  Abbreviate
\begin{center}
$S =$ Span$\lbrace \pi (wz^j) \, | \, w$ is an irreducible word in $F$ and $j \in \Z \rbrace$.
\end{center}
We show $S = A_q(\alpha)$.  Since ${\mathcal F}$ is spanned by its words and since $\pi : {\mathcal F} \rightarrow A_q(\alpha)$ is surjective it suffices to show that $S$ contains the image under $\pi$ of every word in ${\mathcal F}$.  By a {\it counterexample} we mean a word in ${\mathcal F}$ whose image under $\pi$ is not contained in $S$.  We assume there exists a counterexample and obtain a contradiction.  Among all counterexamples let $v$ denote a counterexample with minimal $(x,y)$-length.  Let $t$ denote the $(x,y)$-length of $v$.  Using (\ref{A'q1})--(\ref{A'q3}) we may assume without loss that $v = v'z^{j}$ where $v'$ is a word in $F_t$ and $j \in \Z$.  Recall every word in $F$ of length less than 4 is irreducible.  By construction $v'$ is reducible and so $t \geq 4$.  By Lemma \ref{thm:shaperesult}(iv) there exists $\varpi \in \Omega_t$ and $\lambda \in \Lambda_t$ such that $v' = \varpi + \lambda$.  Now $v = \varpi z^j + \lambda z^j$ so
\begin{eqnarray}
\label{1}
\pi (v) = \pi (\varpi z^j) + \pi (\lambda z^j).
\end{eqnarray}
By construction $\varpi$ is a linear combination of irreducible words so $\pi(\varpi z^{j}) \in S$.  We now show $\pi(\lambda z^{j}) \in S$.  By Lemma \ref{thm:pilambda} and since $\lambda \in \Lambda_t$ we have $\pi (\lambda z^j) \in \pi (F_{t-2} \, z^{j-2})$.  Every word in $F_{t-2} \, z^{j-2}$ has $(x,y)$-length $t-2$ and is therefore not a counterexample by the minimality assumption.  Hence the image under $\pi$ of every word in $F_{t-2} \, z^{j-2}$ is contained in $S$.  Since $F_{t-2} \, z^{j-2}$ is spanned by its words we have $\pi (F_{t-2} \, z^{j-2}) \subseteq S$.  Therefore $\pi (\lambda z^j) \in S$.  We have now shown $\pi(\varpi z^{j}) \in S$ and $\pi(\lambda z^{j}) \in S$ so $\pi (v) \in S$ by (\ref{1}).  This is a contradiction and the result follows.
\hfill $\Box$ \\

\section{A result concerning $(A, \widetilde{A}^*)$-submodules of $V$}

\noindent
Referring to Assumption \ref{assum} and Definition \ref{def:A^*tilde} let $W$ denote an irreducible $(A, \widetilde{A}^*)$-submodule of $V$.  The goal of this section is to prove $V_d \subseteq W$ (see Lemma \ref{thm:VdinW}).  This fact will be used in the proof of Theorem \ref{thm:PVnot0}. \\

\noindent
We note that the arguments given in this section are a modification of the arguments from \cite[Section11]{ItoTer075}.

\begin{definition}
\label{def:endpt}
\rm
With reference to Assumption \ref{assum} and Definition \ref{def:A^*tilde} let $W$ denote an irreducible $(A, \widetilde{A}^*)$-submodule of $V$.  Observe that $W$ is the direct sum of the nonzero spaces among $\lbrace E_iW \rbrace _{i=0}^d$ where $E_i$ is from Definition \ref{def:EiE^*i}(i).
We define
\begin{eqnarray*}
t = \mbox{max} \lbrace i \,|\,0 \leq i \leq d, E_iW \not=0 \rbrace.
\end{eqnarray*}
We call $t$ the {\it endpoint} of  $W$.
\end{definition}

\begin{lemma}
\label{thm:dim1}
With reference to Assumption \ref{assum} and Definition \ref{def:A^*tilde} let $W$ denote an irreducible $(A, \widetilde{A}^*)$-submodule of $V$ and let $t$ denote the endpoint of $W$.  Then \rm{dim}$(E_tW)=1$.
\end{lemma}

\noindent
{\it Proof:} By construction $W$ is an irreducible $(A, \widetilde{A}^*)$-module.  Using this, Lemma \ref{thm:A^*tildediag}, and Lemma \ref{thm:tridiagAA^*tilde} we find that $A|_W, \widetilde{A}^*|_W$ is a $q$-geometric tridiagonal pair on $W$.  Let $s$ denote the diameter of $A|_W, \widetilde{A}^*|_W$.  Note that $\lbrace E_{s-i}W \rbrace _{i=s-t}^{2s-t}$ is a standard ordering of the eigenspaces of $A|_W$.  Applying \cite[Theorem 9.1]{ItoTer075} to $A|_W, \widetilde{A}^*|_W$ we find $\mbox{dim}(E_tW)=1$.
\hfill $\Box $ \\

\begin{lemma}
\label{thm:luzero}
With reference to Assumption \ref{assum} and Definition \ref{def:A^*tilde} let $W$ denote an irreducible $(A, \widetilde{A}^*)$-submodule of $V$ and let $t$ denote the endpoint of $W$.  With reference to Definition \ref{def:EiE^*i}(i) and Definition \ref{def:W_iandfriends} pick $v \in E_tW$ and write $u=G_{d-t}v$.  Then $l u = 0$ where $l$ is the linear transformation from (\ref{rl}).
\end{lemma}

\noindent
{\it Proof:} Observe $u \in W_{d-t}$ by Remark \ref{G_iformulas}.  We assume $d-t \geq 1$; otherwise $l u = 0$ since $l W_0 = 0$.  Observe $l u \in W_{d-t-1}$ by Remark \ref{rlformulas}.  In order to show $l u=0$ we show $l u \in W_{d-t} +\cdots + W_d$.  Using Lemma \ref{thm:AA^*RLrl}(ii) and Definition \ref{def:A^*tilde} we have $c^{-1}l=\widetilde{A}^*-B^{-1}$.  Thus
\begin{eqnarray}
\label{eq:three}
c^{-1}l u = \widetilde{A}^*v - B^{-1}v + c^{-1}l(u-v).
\end{eqnarray}
We are going to show that each of the three terms on the right in (\ref{eq:three}) is contained in $W_{d-t} + \cdots + W_d$.  By the definition of $t$ we have $W = E_0W+ \cdots + E_tW$ so $W \subseteq V_0+ \cdots + V_t$ in view of Definition \ref{def:EiE^*i}(i).  By this and (\ref{AA^*W2}) we find $W \subseteq W_{d-t} + \cdots + W_d$.  By construction $v \in W$ so $\widetilde{A}^*v \in W$.  By these comments $\widetilde{A}^*v \in W_{d-t} + \cdots + W_d$.  We mentioned $v \in W$ so $v \in W_{d-t} + \cdots + W_d$.  Each of $\lbrace W_i \rbrace _{i=d-t}^d$ is an eigenspace for $B^{-1}$ so $B^{-1}v \in W_{d-t} + \cdots + W_d$.  Since $v \in W_{d-t} + \cdots + W_d$ and since $u=G_{d-t}v$ we find $u-v \in W_{d-t+1} + \cdots + W_d$.  Now $c^{-1}l (u-v) \in W_{d-t} + \cdots + W_{d-1}$ so $c^{-1}l (u-v) \in W_{d-t} + \cdots + W_d$.  We have now shown that each of the three terms on the right in (\ref{eq:three}) is contained in $W_{d-t} + \cdots + W_d$.  Therefore $lu \in W_{d-t} + \cdots + W_d$.  Recall $lu \in W_{d-t-1}$.  By this and since $\lbrace W_i \rbrace _{i=0}^d$ is a decomposition of $V$ we find $lu=0$.
\hfill $\Box $ \\

\begin{lemma}
\label{thm:EiW}
With reference to Definition \ref{def:W_iandfriends} and Definition \ref{def:EiE^*i}(i) for $0 \leq i \leq d$ the action of $E_i$ on $W_{d-i}$ coincides with
\begin{eqnarray*}
\sum_{h=0}^{i} \frac{r^{h}} { (q^{2i-d}-q^{2i-d-2}) (q^{2i-d}-q^{2i-d-4}) \cdots (q^{2i-d}-q^{2i-d-2h})}
\end{eqnarray*}
where $r$ is the linear transformation from (\ref{rl}).
\end{lemma}

\noindent
{\it Proof:}  Pick $w \in W_{d-i}$.  We find $E_iw$.  By (\ref{AA^*W2}) and since $E_iw \in V_i$ we find $E_iw \in W_{d-i}+ \cdots +W_d$.  Consequently there exist $w_s \in W_s$ $(d-i \leq s \leq d)$ such that $E_iw=\sum_{s=d-i}^{d}w_s$.  By (\ref{AA^*W1}) and Remark \ref{rlformulas} we have for $0 \leq j \leq d$ that $r$ acts on $W_j$ as $A-q^{d-2j}I$.  Using this and since $(A-q^{2i-d}I)E_i=0$ we find
\begin{eqnarray*}
0 &=& (A-q^{2i-d}I)E_iw \\
  &=& (A-q^{2i-d}I) \sum_{s=d-i}^{d}w_s \\
  &=& \sum_{s=d-i}^{d}(r+q^{d-2s}-q^{2i-d})w_s.
\end{eqnarray*}
Rearranging the terms above we find $0=\sum_{s=d-i+1}^{d}w'_s$ where
\begin{eqnarray*}
w'_s=rw_{s-1} + (q^{d-2s}-q^{2i-d})w_s \qquad (d-i+1 \leq s \leq d).
\end{eqnarray*}
Since $w'_s \in W_s$ for $d-i+1 \leq s \leq d$ and since  $\lbrace W_i \rbrace _{i=0}^d$ is a decomposition of $V$ we find $w'_s=0$ for  $d-i+1 \leq s \leq d$.  Consequently
\begin{eqnarray*}
w_s = (q^{2i-d}-q^{d-2s})^{-1} r w_{s-1} \qquad (d-i+1 \leq s \leq d).
\end{eqnarray*}
By Lemma \ref{thm:bij} and since $w_{d-i}=G_{d-i}E_iw$ we find $w_{d-i}=w$.  From these comments we obtain the desired result.
\hfill $\Box$ \\

\begin{lemma}
\label{thm:updown}
With reference to Assumption \ref{assum} and Definition \ref{def:A^*tilde} let $W$ denote an irreducible $(A, \widetilde{A}^*)$-submodule of $V$ and let $t$ denote the endpoint of $W$.  With reference to Definition \ref{def:EiE^*i}(i) and Definition \ref{def:W_iandfriends} pick $v \in E_tW$ and write $u=G_{d-t}v$.  Then
\begin{eqnarray}
l^i r^i u  \in \mbox{Span}(u) \qquad \qquad (0 \leq i \leq t)
\label{updown1}
\end{eqnarray}
where $r,l$ are the linear transformations from (\ref{rl}).
\end{lemma}

\noindent
{\it Proof:}
We may assume $v \not=0$; otherwise the result is trivial.  Define
\begin{eqnarray}
\label{updown2}
\Delta_i = (\widetilde{A}^*-q^{2t-d}I)(\widetilde{A}^*-q^{2t-d-2}I) \cdots (\widetilde{A}^*-q^{2t-d-2i+2}I).
\end{eqnarray}
Since $\Delta_i$ is a polynomial in $\widetilde{A}^*$ we find $\Delta_iW\subseteq W$. In particular $\Delta_iv \in W$ so $E_t \Delta_iv \in E_tW$.  The vector $v$ spans $E_tW$ by Lemma \ref{thm:dim1} so there exists $m_i \in \K$ such that $E_t \Delta_iv = m_i v$.  By this and since $E_tv=v$ we find $E_t (\Delta_i-m_iI)v=0$.  Now $(\Delta_i-m_iI)v \in E_{0}W + \cdots + E_{t-1}W$ in view of Definition \ref{def:endpt}.  Observe $E_{0}W + \cdots + E_{t-1}W \subseteq V_{0} + \cdots + V_{t-1}$ where the $V_j$ are from Assumption \ref{assum}.  By these comments and (\ref{AA^*W2}) we find $(\Delta_i-m_iI)v \in W_{d-t+1} + \cdots + W_d$.  Consequently $G_{d-t}(\Delta_i-m_iI)v=0$. Recall $G_{d-t}v=u$ so
\begin{eqnarray}
G_{d-t}\Delta_i v=m_iu.
\label{updown3}
\end{eqnarray}
We now evaluate $G_{d-t}\Delta_i v$.  Observe $v=E_tu$ by Lemma \ref{thm:bij} and since $u=G_{d-t} v$.  By Lemma \ref{thm:EiW} there exist nonzero scalars $\gamma_h \in \K$ $(0 \leq h\leq t)$ such that $v=\sum_{h=0}^{t}\gamma_h r^h u$.  For $0 \leq h \leq t$ we compute $G_{d-t}\Delta_i r^h u$.  Keep in mind $r^h u \in W_{d-t+h}$ by Remark \ref{rlformulas}.  First assume $h<i$.  Using Lemma \ref{thm:A^*tildeW} and (\ref{updown2}) we find $\Delta_i r^h u$ is contained in $W_{d-t+h-i}+\cdots+ W_{d-t-1}$ so $G_{d-t}\Delta_i r^h u=0$.  Next assume $h=i$.  Using Lemma \ref{thm:AA^*RLrl}(ii) and Definition \ref{def:A^*tilde} we have $c^{-1}l=\widetilde{A}^*-B^{-1}$ and so $c^{-1}l|_{W_j}=(\widetilde{A}^*-q^{d-2j}I)|_{W_j}$  $(0 \leq j \leq d)$.  Using this and (\ref{updown2}) we find $(\Delta_i-c^{-i}l^i)r^i u$ is contained in $W_{d-t+1}+\cdots + W_{d-t+i}$.  By this and since $c^{-i}l^i r^i u \in W_{d-t}$ we find $G_{d-t}\Delta_i r^i u = c^{-i}l^i r^i u$.  Next assume $h>i$.  Using Lemma \ref{thm:A^*tildeW} and (\ref{updown2}) we find $\Delta_i r^h u$ is contained in $W_{d-t+h-i}+\cdots + W_{d-t+h}$.  By this and since $h>i$ we find $G_{d-t}\Delta_i r^h u =0$.  By these comments we find $G_{d-t}\Delta_i v= \gamma_i c^{-i} l^i r^i u$.  Combining this and (\ref{updown3}) we obtain (\ref{updown1}).
\hfill $\Box $ \\

\begin{lemma}
\label{thm:t=d}
With reference to Assumption \ref{assum} and Definition \ref{def:A^*tilde} let $W$ denote an irreducible $(A, \widetilde{A}^*)$-submodule of $V$ and let $t$ denote the endpoint of $W$.  Then $t=d$.  Moreover, the following holds.  With reference to Definition \ref{def:EiE^*i}(i) and Definition \ref{def:W_iandfriends} pick a nonzero $v \in E_dW$ and write $u=G_0v$.  Let $r,l$ be the linear transformations from (\ref{rl}).  Then $V$ is spanned by the vectors of the form
\begin{eqnarray*}
l^{i_1}r^{i_2}l^{i_3}r^{i_4} \cdots r^{i_n}u
\end{eqnarray*}
where $i_1, i_2, \ldots, i_n$ ranges over all sequences such that $n$ is a nonnegative even integer, and $i_1, i_2, \ldots, i_n$ are integers satisfying $0 \leq i_1 < i_2 <\cdots < i_n \leq d$.
\end{lemma}

\noindent
{\it Proof:}
Let $t$ is the endpoint of $W$.  Pick a nonzero $v \in E_tW$ and write $u=G_{d-t}v$.  Observe $0 \not=u \in W_{d-t}$ by Lemma \ref{thm:bij}.  By Lemma \ref{thm:luzero} and Remark \ref{rlformulas},
\begin{eqnarray}
lu = 0, \qquad \qquad r^{t+1}u=0.
\label{WVdnot01}
\end{eqnarray}
By Lemma \ref{thm:updown},
\begin{eqnarray}
l^i r^i u  \in \mbox{Span}(u) \qquad \qquad (0 \leq i \leq t).
\label{WVdnot02}
\end{eqnarray}
Let $W'$ denote the subspace of $V$ spanned by all vectors of the form
\begin{eqnarray}
l^{i_1}r^{i_2}l^{i_3}r^{i_4} \cdots r^{i_n}u,
\label{WVdnot03}
\end{eqnarray}
where $i_1, i_2, \ldots, i_n$ ranges over all sequences such that $n$ is a nonnegative even integer, and $i_1, i_2, \ldots, i_n$ are integers satisfying $0 \leq i_1 < i_2 <\cdots < i_n\leq t$.  Observe $u \in W'$ so $W' \not=0$.  In order to show $t=d$ we show $W'=V$ and  $W' \subseteq W_{d-t}+\cdots+W_d$.  We now show $W'=V$.  To do this we show that $W'$ is invariant under each of $A,A^*$.  Recall $u \in W_{d-t}$ so $u$ is an eigenvector for $B$ and $B^{-1}$.  Recall the $\K$-algebra $A_q(\alpha)$ from Definition \ref{def:A'qalgebra}.  By Lemma \ref{thm:RLrlBK} and Lemma \ref{thm:RLrlrel} there exists an $A_q(\alpha)$-module structure on $V$ with $\alpha = q^{-4}(q-q^{-1})^3 \lbrack 3 \rbrack !$ where $x,y,z,z^{-1}$ act as $r, l, B, B^{-1}$ respectively.  With respect to this $A_{q}(\alpha)$-module structure we have $W'=A_{q}(\alpha) u$ in view of Lemma \ref{thm:irredform}, Theorem \ref{thm:irredbasisforAq} and (\ref{WVdnot01}), (\ref{WVdnot02}).  It follows that $W'$ is invariant under each of $r, l, B, B^{-1}$.  By Lemma \ref{thm:AA^*RLrl} $A=r+B^{-1}$ and $A^*=l+B+cB^{-1}$.  Using this we have that $W'$ is invariant under each of $A,A^*$.  Since $A,A^*$ is a tridiagonal pair on $V$ we have $W'=V$.  We now show $W' \subseteq W_{d-t} + \cdots + W_d$.  By Remark \ref{rlformulas} the vector (\ref{WVdnot03}) is contained in $W_{d-t+i}$ where $i=\sum_{h=1}^n i_h(-1)^h$.  From the construction $0\leq i \leq t$ so $W_{d-t+i} \subseteq W_{d-t}+\cdots + W_d$.  Therefore the vector (\ref{WVdnot03}) is contained in $W_{d-t}+\cdots + W_d$ so $W'\subseteq W_{d-t} +\cdots + W_d$.  We have shown $W'=V$ and $W' \subseteq W_{d-t} + \cdots + W_d$.  Since  $\lbrace W_i \rbrace _{i=0}^d$ is a decomposition of $V$ we find $t=d$ and the result follows.
\hfill $\Box $ \\

\begin{lemma}
\label{thm:VdinW}
With reference to Assumption \ref{assum} and Definition \ref{def:A^*tilde} let $W$ denote an irreducible $(A, \widetilde{A}^*)$-submodule of $V$.  Then $V_d \subseteq W$.
\end{lemma}

\noindent
{\it Proof:} Recall by Lemma \ref{thm:t=d} that the endpoint of $W$ is $d$.  So $E_dW \not=0$ by Definition \ref{def:endpt}.  We have $U_0=W_0$ by (\ref{AA^*U2}),(\ref{AA^*W2}) and so $\mbox{dim}(W_0)=1$ by Lemma \ref{thm:dim1U0}.  Using this and Lemma \ref{thm:bij} we find $\mbox{dim}(V_d)=1$.  We have $0 \neq E_dW \subseteq V_d$ so $E_dW=V_d$.  But $E_dW \subseteq W$ by (\ref{EE^*AA^*1}) so $V_d \subseteq W$.
\hfill $\Box$ \\

\section{$A, \widetilde{A}^* $ is a tridiagonal pair}

\noindent
In this section we show $A, \widetilde{A}^*$ is a $q$-geometric tridiagonal pair of $V$ if and only if \\ $P( q^{2d-2} (q-q^{-1})^{-2} ) \not=0$.  The proof of this depends on the following lemma.

\begin{lemma}
\label{thm:PVnot0}
With reference to Assumption \ref{assum}, Definition \ref{def:PV}, and Definition \ref{def:A^*tilde} the following are equivalent:
\begin{enumerate}
\item[{\rm (i)}]
$V$ is irreducible as an $(A, \widetilde{A}^*)$-module.
\item[{\rm (ii)}]
$P( q^{2d-2} (q-q^{-1})^{-2} ) \neq 0$.
\end{enumerate}
\end{lemma}

\noindent
{\it Proof:} (i) $\Rightarrow$ (ii) We assume $P( q^{2d-2} (q-q^{-1})^{-2} )=0$ and derive a contradiction.  Define
\begin{eqnarray*}
X_i=(V_i+ \cdots + V_d) \cap (\widetilde{V}^*_{d-i+1}+ \cdots + \widetilde{V}^*_d) \qquad (1 \leq i \leq d)
\end{eqnarray*}
where the $V_j$ are from Assumption \ref{assum} and the $\widetilde{V}^*_j$ are from Definition \ref{def:V^*_itilde}.  Further define $X=X_1+ \cdots + X_d$.  We will show that $X$ is an $(A, \widetilde{A}^*)$-submodule of $V$ and $X \neq V$, $X \neq 0$.  We first show $AX \subseteq X$.  For $1 \leq i \leq d$ we have $(A-q^{2i-d}I) \sum_{j=i}^d V_j = \sum_{j=i+1}^d V_j$ by Assumption \ref{assum} and  $(A-q^{2i-d}I) \sum_{j=d-i+1}^d \widetilde{V}^*_j \subseteq \sum_{j=d-i}^d \widetilde{V}^*_j$ by Lemma \ref{thm:tridiagAA^*tilde}(ii).  By these comments
\begin{eqnarray*}
(A-q^{2i-d}I)X_i \subseteq X_{i+1} \quad (1 \leq i \leq d-1), \quad (A-q^dI)X_d=0
\end{eqnarray*}
and it follows $AX \subseteq X$.  We now show $\widetilde{A}^*X \subseteq X$.  For $1 \leq i \leq d$ we have \\ $(\widetilde{A}^*-q^{2i-d-2}I) \sum_{j=i}^d V_j = \sum_{j=i-1}^d V_j$ by Lemma \ref{thm:tridiagAA^*tilde}(i) and \\ $(\widetilde{A}^*-q^{2i-d-2}I) \sum_{j=d-i+1}^d \widetilde{V}^*_j \subseteq \sum_{j=d-i+2}^d \widetilde{V}^*_j$ by Definition \ref{def:V^*_itilde}.  By these comments
\begin{eqnarray*}
(\widetilde{A}^*-q^{2i-d-2}I)X_i \subseteq X_{i-1} \quad (2 \leq i \leq d-1), \quad (\widetilde{A}^*-q^{-d}I)X_1=0
\end{eqnarray*}
and it follows $\widetilde{A}^*X \subseteq X$.  We have now shown that $X$ is an $(A, \widetilde{A}^*)$-submodule of $V$.  We now show $X \neq V$.  For $1 \leq i \leq d$ we have $X_i \subseteq V_i + \cdots + V_d$ so $X_i \subseteq V_1 + \cdots + V_d$.  It follows $X \subseteq V_1 + \cdots +V_d$ and so $X \neq V$.  We now show $X \neq 0$.  To do this we display a nonzero vector in $X_d$.  Pick a nonzero vector $u \in U_0$.  Applying Theorem \ref{thm:PVformula} we find $\widetilde{E}^*_0E_du=0$.  Write $v=E_du$ and notice $v \in V_d$.  By (\ref{AA^*U2}) and (\ref{AA^*W2}) we find $U_0=W_0$ and so $v \neq 0$ by Lemma \ref{thm:bij}.  Observe $\widetilde{E}^*_0v=0$ so $v \in \widetilde{V}^*_1+ \cdots + \widetilde{V}^*_d$ by Definition \ref{def:EiE^*i}(ii).  From these comments $v \in X_d$.  We have displayed a nonzero vector $v$ contained in $X_d$.  Of course $X_d \subseteq X$ so $X \neq 0$.  We have now shown that $X$ is an $(A, \widetilde{A}^*)$-submodule of $V$ and $X \neq V$, $X \neq 0$.  This contradicts our assumption that $V$ is irreducible as an $(A, \widetilde{A}^*)$-module.  We conclude $P( q^{2d-2} (q-q^{-1})^{-2} ) \not=0$.  \\
(ii) $\Rightarrow$ (i) Let $W$ denote an irreducible $(A, \widetilde{A}^*)$-submodule of $V$.  We show $W=V$.  To do this we show $W$ is invariant under each of $A,A^*$.  By construction $W$ is invariant under $A$.  In order to show $W$ is invariant under $A^*$ we show $W$ is invariant under $B$.  We define $\overline {W} = \{ \, w \in W \, | \, Bw \in W \, \}$ and show $\overline {W}=W$.  Using (\ref{AA^*B1}) we find $A \overline {W} \subseteq \overline {W}$.  Using (\ref{A^*tildeBrel1}) we find $\widetilde{A}^* \overline {W} \subseteq \overline {W}$.  We now show $\overline {W} \neq 0$.  Let $0 \neq v \in V_d$.  By Lemma \ref{thm:VdinW} $v\in W$.  Combining (\ref{AA^*U2}) and (\ref{AA^*W2}) we find $U_0=W_0$.  Using this and Lemma \ref{thm:bij} (with $i=d$) we have that $E_d: U_0 \rightarrow V_d$ is a bijection.  So there exists a nonzero $u \in U_0$ such that $E_du=v$.  Using Theorem \ref{thm:PVformula} we find $\widetilde{E}^*_0v=q^{2d(1-d)} c^{-d} P( q^{2d-2} (q-q^{-1})^{-2} )u$.  Since $P( q^{2d-2} (q-q^{-1})^{-2} )$ and $u$ are both nonzero we find $\widetilde{E}^*_0v \neq 0$.  Using (\ref{EE^*AA^*2}) and since $v \in W$ we have $\widetilde{E}^*_0v \in W$.  Using Lemma \ref{thm:A^*tildeW} and Lemma \ref{thm:A^*tildediag} we find $\widetilde{V}^*_0=W_0$.  Hence $\widetilde{E}^*_0v \in W_0$ and so $\widetilde{E}^*_0v \in \overline {W}$.  By these comments we find $\overline {W} \neq 0$.  We have now shown $\overline {W}$ is nonzero and invariant under each of $A,\widetilde{A}^*$.  Therefore $\overline{W}=W$ since $W$ is an irreducible $(A, \widetilde{A}^*)$-module.  We have now shown $W$ is invariant under $B$.  By construction $W$ is invariant under $\widetilde{A}^*$.  So by Definition \ref{def:A^*tilde} $W$ is invariant under $A^*$.  We now know that $W$ is nonzero and invariant under each of $A,A^*$.  Since $A,A^*$ is a tridiagonal pair on $V$ we find $W=V$ and the result follows.
\hfill $\Box$ \\

\begin{lemma}
\label{thm:AA^*tildetrip}
With reference to Assumption \ref{assum}, Definition \ref{def:PV}, and Definition \ref{def:A^*tilde} the following holds.  $A, \widetilde{A}^*$ is a $q$-geometric tridiagonal pair on $V$ if and only if $P( q^{2d-2} (q-q^{-1})^{-2} ) \not=0$.
\end{lemma}

\noindent
{\it Proof:}  Immediate from Assumption \ref{assum}, Lemma \ref{thm:A^*tildediag}, Lemma \ref{thm:tridiagAA^*tilde}, and Lemma \ref{thm:PVnot0}.
\hfill $\Box$ \\

\section{The proof of Theorem \ref{thm:main}}

\noindent
In this section we give a proof of our main result by providing the required action of $\boxtimes_q$ on $V$. \\

\noindent
{\it Proof of Theorem \ref{thm:main}:}  (i) $\Rightarrow$ (ii)  By \cite[Theorem 12.3]{ItoTer072} the action of $x_{30}$ on $V$ is diagonalizable and the set of distinct eigenvalues is $\{ q^{2i-d} \, | \, 0 \leq i \leq d \, \}$.  For $0 \leq i \leq d$ let $X_i$ denote the eigenspace of $x_{30}$ corresponding to the eigenvalue $q^{2i-d}$.  First we show $x_{30}$ acts on $V$ as $B$ from Definition \ref{def:B}.  By Definition \ref{def:qtet} $qx_{01}x_{30}-q^{-1}x_{30}x_{01}=(q-q^{-1})I$ and so we have $(x_{01}-q^{d-2i}I)X_i \subseteq X_{i+1}$ for $0 \leq i \leq d$ in view of Lemma \ref{thm:bidiag2}.  Using Definition \ref{def:qtet} we find $qx_{30}(x_{30}+cx_{23})-q^{-1}(x_{30}+cx_{23})x_{30}=(q-q^{-1})(x^2_{30}+cI)$ and so $(x_{30}+cx_{23}-q^{2i-d}I-cq^{d-2i}I)X_i \subseteq X_{i-1}$ for $0 \leq i \leq d$ in view of Lemma \ref{thm:bidiag3}.  Recall by construction $x_{01}$ acts as $A$ and $x_{30}+cx_{23}$ acts as $A^*$.  Using these comments and \cite[Theorem~4.6]{ItoTanTer01} (with $V_i$ replaced by $V_{d-i}$) we find  $X_i=(V^*_0 + \cdots + V^*_i) \cap (V_0 + \cdots + V_{d-i})$ for $0 \leq i \leq d$.  Using this and Remark \ref{W_iformulas} we find $X_i=W_i$ for $0 \leq i \leq d$.  In view of Definition \ref{def:B} we have now shown $x_{30}$ acts on $V$ as $B$.  Using this and Definition \ref{def:A^*tilde} we find $x_{23}$ acts as $\widetilde{A}^*$.  By \cite[Theorem 10.3]{ItoTer072} and \cite[Theorem 2.7]{ItoTer075} we find $A, \widetilde{A}^*$ is a $q$-geometric tridiagonal pair of $V$.  So $P( q^{2d-2} (q-q^{-1})^{-2} ) \not=0$ by Lemma \ref{thm:AA^*tildetrip}. \\
(ii) $\Rightarrow$ (i) By Lemma \ref{thm:AA^*tildetrip} $A, \widetilde{A}^*$ is a $q$-geometric tridiagonal pair on $V$.  Using this, \cite[Theorem 2.7]{ItoTer075}, and \cite[Theorem 10.4]{ItoTer072} there exists a unique irreducible $\boxtimes_q$-module structure on $V$ such that $x_{01}$ acts as $A$ and $x_{23}$ acts as $\widetilde{A}^*$.  It remains to show that $x_{30}+cx_{23}$ acts as $A^*$.  To do this we show $x_{30}$ acts on $V$ as $B$ from Definition \ref{def:B}.  By \cite[Theorem 12.3]{ItoTer072} the action of $x_{30}$ on $V$ is diagonalizable and the set of distinct eigenvalues is $\{ q^{2i-d} \, | \, 0 \leq i \leq d \, \}$.  For $0 \leq i \leq d$ let $X_i$ denote the eigenspace of $x_{30}$ corresponding to the eigenvalue $q^{2i-d}$.  Using \cite[Theorem 16.4]{ItoTer072} we find $X_i=(\widetilde{V}^*_0 + \cdots + \widetilde{V}^*_i) \cap (V_0 + \cdots + V_{d-i})$ for $0 \leq i \leq d$.  Recall by Definition \ref{def:B} that for $0 \leq i \leq d$, $W_i$ is the eigenspace of $B$ with eigenvalue $q^{2i-d}$.  We show $X_i=W_i$ for $0 \leq i \leq d$.  Combining Lemma \ref{thm:qweylsum} and (\ref{A^*tildeBrel1}) we find $W_0+ \cdots + W_i=\widetilde{V}^*_0+ \cdots + \widetilde{V}^*_i$ for $0 \leq i \leq d$.  Using this and (\ref{AA^*W2}) we find $V^*_0+ \cdots + V^*_i=\widetilde{V}^*_0 + \cdots + \widetilde{V}^*_i$ for $0 \leq i \leq d$.  Using this and Remark \ref{W_iformulas} we find $X_i=W_i$ for $0 \leq i \leq d$.  We have now shown $x_{30}$ acts on $V$ as $B$.  Recall by construction $x_{23}$ acts as $\widetilde{A}^*$.  Therefore by Definition \ref{def:A^*tilde} $x_{30}+cx_{23}$ acts as $A^*$ and the result follows.
\hfill $\Box$ \\

\noindent
For the sake of completeness we now make a few comments regarding the $\boxtimes_q$-module structure on $V$ given in Theorem \ref{thm:main}.

\begin{lemma}
\label{thm:boxqaction}
With reference to Assumption \ref{assum}, Theorem \ref{thm:main}, and Definition \ref{def:V^*_itilde} the following holds.  For each generator $x_{ij}$ of $\boxtimes_q$ and for $0 \leq n \leq d$ the eigenspace of $x_{ij}$ corresponding to the eigenvalue $q^{2n-d}$ is given as follows.
\medskip
\begin{center}
\begin{tabular}[t]{c|c}
       {\rm generator} & {\rm eigenspace corresponding to eigenvalue $q^{2n-d}$}
 \\ \hline  \hline
	$x_{01}$ & $V_n$  \\
	$x_{23}$ & $\widetilde{V}^*_{d-n}$   \\
        $x_{30}$ & $(\widetilde{V}^*_0+ \cdots +\widetilde{V}^*_n) \cap (V_0+ \cdots +V_{d-n})$ \\
        $x_{12}$ & $(\widetilde{V}^*_n+ \cdots +\widetilde{V}^*_d) \cap (V_{d-n}+ \cdots +V_d)$ \\
        $x_{31}$ & $(\widetilde{V}^*_0+ \cdots +\widetilde{V}^*_n) \cap (V_n+ \cdots +V_d)$  \\
        $x_{13}$ & $(\widetilde{V}^*_0+ \cdots +\widetilde{V}^*_{d-n}) \cap (V_{d-n}+ \cdots +V_d)$  \\
        $x_{20}$ & $(\widetilde{V}^*_{d-n}+ \cdots +\widetilde{V}^*_d) \cap (V_0+ \cdots +V_{d-n})$ \\
	$x_{02}$ & $(\widetilde{V}^*_n+ \cdots +\widetilde{V}^*_d) \cap (V_0+ \cdots +V_n)$
        \end{tabular}
\end{center}
\end{lemma}

\noindent
{\it Proof:} The first row of the table follows immediately from Theorem \ref{thm:main}.  In the proof of Theorem \ref{thm:main} we showed $x_{23}$ acts on $V$ as $\widetilde{A}^*$.  Using this we immediately obtain the second row of the table.  The remaining six rows of the table follow from the first two rows and \cite[Theorem 16.4]{ItoTer072}.
\hfill $\Box$ \\

\begin{lemma}
With reference to Assumption \ref{assum} and the $\boxtimes_q$-module structure on $V$ given in Theorem \ref{thm:main} the following holds.
\begin{enumerate}
\item $x_{30}$ acts on $V$ as $B$.
\item $x_{23}$ acts on $V$ as $\widetilde{A}^*$.
\item $x_{31}$ acts on $V$ as $K$.
\item $x_{13}$ acts on $V$ as $K^{-1}$.
\end{enumerate}
where $B,\widetilde{A}^*,K$ are from Definition \ref{def:B}, Definition \ref{def:A^*tilde}, and Definition \ref{def:K} respectively.
\end{lemma}

\noindent
{\it Proof:} (i),(ii):  These were shown in the proof of Theorem \ref{thm:main}. \\
(iii),(iv): Recall for $0 \leq i \leq d$ that $U_i$ is the eigenspace for $K$ corresponding to the eigenvalue $q^{2i-d}$.  In the proof of Theorem \ref{thm:main} we showed for $0 \leq i \leq d$ that $V^*_0+ \cdots + V^*_i=\widetilde{V}^*_0 + \cdots + \widetilde{V}^*_i$.  Using this and Lemma \ref{def:U_i} we find $U_i=(\widetilde{V}^*_0 + \cdots + \widetilde{V}^*_i) \cap (V_i + \cdots + V_d)$ for $0 \leq i \leq d$.  The result now follows from rows five and six of the table in Lemma \ref{thm:boxqaction}.
\hfill $\Box$ \\

\section{Acknowledgment}

\noindent
The research which lead to this paper was done at the University of Wisconsin -- Madison as part of the author's PhD thesis.  This thesis was written under the direction of Paul Terwilliger.  The author would like to express his gratitude to Paul Terwilliger for his many helpful suggestions.

\bibliographystyle{amsplain}
\bibliography{mybib}

\noindent
Darren Funk-Neubauer \hfil\break
\noindent Department of Mathematics and Physics \hfil\break
\noindent Colorado State University - Pueblo \hfil\break
\noindent 2200 Bonforte Boulevard \hfil\break
\noindent Pueblo, CO 81001 USA \hfil\break
email:   {\tt darren.funkneubauer@colostate-pueblo.edu} \hfil\break
phone:  (719) 549 - 2693 \hfil\break
fax:  (719) 549 - 2962

\end{document}